\documentclass{article}

\usepackage{multicol}        
\usepackage[bottom]{footmisc}
\usepackage{amssymb}
\usepackage{amsmath}
\usepackage{mathrsfs}
\usepackage{color}
\usepackage{caption}
\usepackage[normalem]{ulem}
\usepackage{calligra}
\usepackage[T1]{fontenc}
\usepackage{lipsum}
\usepackage{xcolor}

\newtheorem{theorem}{Theorem}
\newtheorem{lemma}[theorem]{Lemma}

\newcommand{\be}{\begin{equation}}
\newcommand{\ee}{\end{equation}}
\newcommand{\bea}{\begin{eqnarray}}
\newcommand{\eea}{\end{eqnarray}}
\newcommand{\bean}{\begin{eqnarray*}}
\newcommand{\eean}{\end{eqnarray*}}
\newcommand{\la}{\label}
\newcommand{\ba}{\begin{array}}
\newcommand{\ea}{\end{array}}

\newcommand\blfootnote[1]{%
  \begingroup
  \renewcommand\thefootnote{}\footnote{#1}%
  \addtocounter{footnote}{-1}%
  \endgroup
}

\newcommand{\xO}{\Omega}

\newcommand{\R}{{\mathbb{R}}}

\newcommand{\Hp}{{\mathbb{H}}}
 
 \newcommand{\cC}{{\cal C}}

 \newcommand{\cD}{{\cal D}}

 \newcommand{\as}{\mbox{ as }}

 \newcommand{\inprod}[2]{{\langle{#1},{#2}\rangle}}

 \newcommand{\darr}[4]{{\left\{\begin{array}{ll}
   {#1}&{#2}\\[0.2cm]
   {#3}&{#4}
 \end{array}\right.}}

\newcommand{\ia}{({\rm i})}
\newcommand{\ib}{({\rm ii})}

\def\Int{\displaystyle\int}

\def\Frac{\displaystyle\frac}

\newcommand{\cic}{{C^{\infty}_c}}

\newcommand{\sph}{{  {\rm S}^{n-1}}}
\newcommand{\sphp}{{  {\rm S}^{n-1}_+}}

\newcommand{\us}{ {\rm S}^{n-1}}

\title{Best Sobolev constants in the presence of sharp Hardy terms in Euclidean and hyperbolic space}
\date{}

\author{
G. Barbatis\thanks{Department of Mathematics,
 National and Kapodistrian University of Athens,  15784 Athens, Greece}
 \and A. Tertikas
\thanks{Department of Mathematics and Applied Mathematics,
 University of Crete, 70013 Heraklion, Greece }
}

\begin{document}

\date{\today}

\maketitle
\blfootnote{Email addresses: gbarbatis@math.uoa.gr; tertikas@uoc.gr}


\begin{abstract}
\noindent

In this article we compute the best Sobolev constants for various Hardy-Sobolev inequalities with sharp Hardy term.
This is carried out in three different environments: interior point singularity in Euclidean space, interior point singularity in hyperbolic space and boundary point singularity in Euclidean domains.

\end{abstract}

\vspace{11pt}


\vspace{6pt}

\section{Introduction}

The standard Hardy inequality in $\R^n$, $n\geq 3$, reads
\be
\int_{\R^n}|\nabla u|^2dx \geq\Big( \frac{n-2}{2}\Big)^2 \int_{\R^n}\frac{u^2}{|x|^2}dx \;  , \qquad u\in \cic(\R^n).
\la{eq:hi}
\ee
The constant $((n-2)/2)^2$ is sharp and is not attained in any reasonable function space such as $\cD^{1,2}(\R^n)$, the completion of $\cic(\R^n)$ with respect to the norm $\|\nabla u\|_{L^2}$.
The same remains true if we replace $\R^n$ by $B_1$, the unit ball centered at zero.

Similarly the Sobolev inequality in $\R^n$, $n\geq 3$, reads
\be
\int_{\R^n}|\nabla u|^2dx 
\geq S_n  \Big( \int_{\R^n}|u|^{2^*} dx \Big)^{2/2^*} \; , \qquad  u\in \cic(\R^n),
\la{sobolev}
\ee
where $2^* =2n/(n-2)$ and the sharp constant $S_n$ is given by
\[
S_n = \pi n (n-2) \Big( \frac{\Gamma(\frac{n}{2})}{\Gamma(n)}   \Big)^{\frac{2}{n}} .
\]
This inequality has as a minimizer in $\cD^{1,2}(\R^n)$ the function
\[
u(x) =  \big( 1 +|x|^2 \big)^{-\frac{n-2}{2}} \; , \quad x\in\R^n,
\]
as well as translates and scaled versions of it.

The following family of inequalities interpolate between the Hardy inequality (\ref{eq:hi}) and the Sobolev inequality:
for any $2 < p \leq 2^*$ there holds
\be
\int_{\R^n}|\nabla u|^2 dx \geq S_{n,p} \Big(  \int_{\R^n}|x|^{ \frac{p(n-2)}{2}-n}|u|^p dx\Big)^{2/p} , \qquad u\in C^{\infty}_c(\R^n).
\la{ws}
\ee
The sharp constant $S_{n,p}$ has been computed in \cite{Lieb} and is given by
\be
S_{n,p}=  2p \Big( \frac{n-2}{2}  \Big)^{\frac{p+2}{2}}  
 \bigg[ 
\frac{ 2 \pi^{n/2}\Gamma^2(\frac{p}{p-2})}{ (p-2) \Gamma(\frac{n}{2}) \Gamma(\frac{2p}{p-2}) } 
 \bigg]^{\frac{p-2}{p}}  \; \; , \qquad 2<p \leq 2^* \; ,
\la{lieb}
\ee
and one minimizer is the function
\be
u(x) = \big(  1+|x|^{\frac{(p-2)(n-2)}{2}} \big)^{-\frac{2}{p-2}} , \qquad x\in\R^n .
\la{minim}
\ee

Let us define
\[
X( t) := \frac{1}{1 -\ln t}   \; , \quad t\in (0,1).
\]
In \cite{AFT}, following earlier work in \cite{FT}, the Hardy-Sobolev inequality
\be
\int_{B_1}|\nabla u|^2dx -\Big( \frac{n-2}{2}\Big)^2 \int_{B_1}\frac{u^2}{|x|^2}dx
\geq
(n-2)^{-\frac{2(n-1)}{n}}S_n
 \Big( \int_{B_1} X^{\frac{2(n-1)}{n-2}}(\alpha |x|) |u|^{2^*} dx \Big)^{2/2^*} 
\la{intro:1}
\ee
was established for all $0<\alpha\leq \alpha_n$ and all $u\in \cic(B_1)$, where
$\alpha_n= e^{\frac{n-3}{n-2}}$. The exponent $2(n-1)/(n-2)$ is sharp which in particular implies the necessity of the logarithmic factor $X$
for the validity of (\ref{intro:1}).
Moreover it was shown that the Sobolev constant is sharp for any $0<\alpha\leq \alpha_n$.

In this work we prove the following sharp interpolated inequality: for any
$2<p \leq 2^*$, $0<\alpha \leq \alpha_n$ and all $u\in\cic(B_1)$ there holds
\be
\int_{B_1}|\nabla u|^2dx -\Big( \frac{n-2}{2}\Big)^2 \int_{B_1}\frac{u^2}{|x|^2}dx
\geq  (n-2 )^{-\frac{p+2}{p}}S_{n,p}
 \Big( \int_{B_1} |x|^{\frac{p(n-2)}{2} -n}  X(\alpha |x|)^{\frac{p+2}{2}} |u|^{p}  dx \Big)^{2/p} .
\la{intro:mel1}
\ee
Actually in Theorems \ref{thm:2019} and \ref{thm7} we establish improved versions of (\ref{intro:mel1}).

More generally, our aim in this work is to obtain inequalities analogous to (\ref{intro:mel1})
in different geometric contexts, namely Euclidean or hyperbolic with interior point singularity and Euclidean with boundary point singularity.

One geometric environment where there has been a lot of recent activity on Hardy and Sobolev inequalities is the hyperbolic space $\Hp^n$; 
see \cite{AK,BGG,BGGP,CGMS,GMR,H,KO,LY,MS,Ng,ST}.
The analogue of the Sobolev inequality (\ref{sobolev}) in the hyperbolic space reads 
\[
\int_{\Hp^n}|\nabla_{\Hp^n} u|^2dV - \frac{n(n-2)}{4} \int_{\Hp^n}u^2 dV 
\geq
S_{n}\Big( \int_{\Hp^n}  |u|^{2^*} dV \Big)^{2/2^*},
\]
and the constant $S_n$ is sharp \cite{H}. In fact the full interpolation inequality (\ref{ws}) can be translated to the hyperbolic environment where it takes the form
\[
\int_{\Hp^n}|\nabla_{\Hp^n} u|^2dV - \frac{n(n-2)}{4} \int_{\Hp^n}u^2 dV 
\geq
S_{n,p}\Big( \int_{\Hp^n} (\sinh \rho)^{\frac{p(n-2)}{2}-n} |u|^{p} dV \Big)^{2/p},
\]
where $S_{n,p}$ is given by (\ref{lieb}) and is sharp; see also \cite[Corollary 2.3]{KO} for an analogous interpolated result.

It is not difficult to see that inequality (\ref{intro:mel1}) can be transformed to the hyperbolic environment giving that
for all $2<p\leq 2^*$, $0<\alpha\leq \alpha_n$ and all $v\in \cic(\Hp^n)$ there holds
\[
\int_{\Hp^n}|\nabla_{\Hp^n} v|^2dV -\Big( \frac{n-2}{2}\Big)^2 \int_{\Hp^n}\frac{v^2}{\rho^2}dV 
\geq
(n-2)^{-\frac{p+2}{p}}S_{n,p}\Big( \int_{\Hp^n}  (\sinh \rho)^{\frac{p(n-2)}{2}-n} X^{\frac{p+2}{2}}(\alpha \tanh (\rho/2)) |v|^{p} dV \Big)^{2/p}
\]
Moreover the Sobolev constant is sharp for any $0<\alpha\leq\alpha_n$.

Actually we prove a slightly stronger result that reads
\begin{theorem} {\bf (Hardy-Sobolev inequality)}
Let $n\geq 3$ and $2<p\leq 2^*$. For any $0<\alpha\leq\alpha_n$ there holds
\bea
&&\hspace{-1cm}\int_{\Hp^n}|\nabla_{\Hp^n} u|^2dV -\frac{n(n-2)}{4} \int_{\Hp^n}u^2dV   -\Big( \frac{n-2}{2}\Big)^2 \int_{\Hp^n}\frac{u^2}{\sinh ^2\rho}dV
 \nonumber \\[0.2cm]
&& \hspace{2.5cm}\geq
(n-2)^{-\frac{p+2}{p}}S_{n,p}\Big( \int_{\Hp^n}  (\sinh \rho)^{\frac{p(n-2)}{2}-n} X^{\frac{p+2}{2}}(\alpha \tanh (\rho/2)) |u|^{p} dV \Big)^{2/p},
\label{101}
\eea
for all $u\in\cic(\Hp^n)$. Moreover the constant $(n-2)^{-\frac{p+2}{p}}S_{n,p}$ is sharp.
\label{thm:intro1}
\end{theorem}
It is worth noting that in case $n=3$ the Sobolev constant of inequality (\ref{101})  is equal to
$S_{3,p}$ whereas for $n\geq 4$ it is strictly smaller than $S_{n,p}$.

In a different direction in \cite{AK,BGG} the following non-improvable \cite{BGG} Poincar\'{e}-Hardy inequality was established
\[
\int_{\Hp^n}|\nabla_{\Hp^n} u|^2dV - \Big(\frac{n-1}{2}\Big)^2 \int_{\Hp^n}u^2 dV 
\geq \frac{1}{4}  \int_{\Hp^n}\frac{u^2}{\rho^2} dV + \frac{(n-1)(n-3)}{4} \int_{\Hp^n}\frac{u^2}{\sinh^2\rho} dV.
\]
We note that $(n-1)^2/4$ is the bottom of the spectrum of the Laplace operator on $\Hp^n$.
Here we consider the Poincar\'{e}-Hardy inequality
\be
 \int_{\Hp^n}|\nabla_{\Hp^n} v|^2dV  \geq \Big( \frac{n-1}{2}\Big)^2 \int_{\Hp^n} v^2 dV 
+  \Big( \frac{n-2}{2}\Big)^2 \int_{\Hp^n} \frac{v^2}{ \sinh^2\rho}  dV \, ,
\la{889}
\ee
as well as, for $n\geq 3$, the Poincar\'{e}-Sobolev inequality
\be
 \int_{\Hp^n}|\nabla_{\Hp^n} v|^2dV  \geq \Big( \frac{n-1}{2}\Big)^2 \int_{\Hp^n} v^2 dV 
+ \overline{S}_{n,p}\Big( \int_{\Hp^n}(\sinh \rho)^{\frac{p(n-2)}{2}-n} |v|^{p}  dV \Big)^{2/p} , \quad\quad v\in \cic(\Hp^n),
\la{intro:100}
\ee
where $2<p\leq 2^*$ and $\overline{S}_{n,p}$ denotes the best constant. The positivity of $\overline{S}_{n,p}$ follows from the positivity
of $\overline{S}_{n,2^*}$  (see \cite{MS}) together with (\ref{889}).

Actually, using the half-space or the unit ball model of $\Hp^n$ one can see that inequality
(\ref{intro:100}) can be written in two equivalent ways, namely
\be
 \int_{\R^n_+}|\nabla u|^2dx  \geq \frac{1}{4}  \int_{\R^n_+} \frac{u^2}{x_n^2}dx 
+ \overline{S}_{n,p}\bigg( \int_{\R^n_+}  \Big(  \frac{ |x-e_n| \; |x+e_n|}{2} \Big)^{ \frac{(n-2)p}{2}-n} |u|^{p}  dx \bigg)^{2/p} , \quad\quad u\in \cic(\R^n_+),
\la{intro:101}
\ee
and
\be
 \int_{B_1}|\nabla u|^2dx  \geq \frac{1}{4} \int_{B_1} \frac{u^2}{ \big( \frac{1-|x|^2}{2}  \big)^2 }  dx 
+ \overline{S}_{n,p}\Big( \int_{B_1} |x|^{ \frac{p(n-2)}{2
}-n}|u|^{p}  dx \Big)^{2/p} , \quad\quad u\in \cic(B_1),
\la{intro:210}
\ee
the best constants of (\ref{intro:100}), (\ref{intro:101}) and (\ref{intro:210}) being equal. 
For work related to inequality (\ref{intro:101}) see also \cite{FMT,TT}.

The precise value of $\overline{S}_{n,p}$ is not known in general. In the case $n=3$, $p=2^*=6$ it has been shown in \cite{BFL} that
\[
\overline{S}_{3,6} =S_3 =3 \big(\frac{\pi}{2} \big)^{4/3}.
\]
Adapting the ideas of \cite{BFL} we compute in Theorem \ref{thm15} the constant $\overline{S}_{3,p}$ for any $2<p < 6$ and find that
\be
\overline{S}_{3,p} =S_{3,p} = \frac{p}{2^{\frac{2}{p}}} 
\bigg[ 
\frac{ 4 \pi \Gamma^2(\frac{p}{p-2})}{ (p-2)  \Gamma(\frac{2p}{p-2}) } 
 \bigg]^{\frac{p-2}{p}}.
\la{identify}
\ee

The next problem we address is how the Sobolev constant in (\ref{intro:100}) is affected when we add a
Hardy term with sharp constant in the RHS. To answer this question we need to study in detail the existence and asymptotic behaviour of positive solutions of the following two problems:
\be
\darr{g''(t) + \Frac{1}{4\sinh^2 t}g(t)=0 ,}{t>0,}{ \lim_{t \to +\infty}g(t) =1,}{}
\la{intro:g}
\ee
and, for $n\geq 3$,
\[
\darr{h''(t) - \Frac{(n-1)(n-3)}{4\sinh^2 t}h(t)=0,}{t>0,}{  \lim_{t \to +\infty}h(t) =1 \, .}{}
\]
We shall see that these problems have unique solutions $g$ and $h$ which are actually positive and behave near zero in a way that
allows us to define a function $\rho=\rho(t)$ by
\be
 \int_0^{\rho(t)} \frac{dr}{h(r)^2} = \int_0^{t} \frac{ds}{g(s)^2}  \;  , \qquad t>0.
\la{intro:fnr}
\ee
We then define
\be
Y(t) = (n-2)\frac{h\big(\rho(t)\big)^2 \sinh t }{g(t)^2  \sinh \rho(t) } \; , \qquad t>0.
\la{def:Y}
\ee
We have the following
\begin{theorem}
{\bf (Poincar\'{e}-Hardy-Sobolev inequality I)} Let $n\geq 3$ and $2<p\leq 2^*$.  There holds
\bea
\int_{\Hp^n}|\nabla_{\Hp^n} v|^2dV &\geq& \Big( \frac{n-1}{2}\Big)^2 \int_{\Hp^n}v^2 dV +
\Big( \frac{n-2}{2}\Big)^2 \int_{\Hp^n}\frac{v^2}{\sinh^2\rho}dV  \nonumber \\[0.2cm]
&& + (n-2)^{-\frac{p+2}{p}} \overline{S}_{n,p}\Big( \int_{\Hp^n} (\sinh\rho)^{\frac{p(n-2)}{2}-n} 
Y^{\frac{p+2}{2}}(\rho)|v|^{p} dV \Big)^{2/p} ,
\la{intro:102}
\eea
for all $ v\in \cic(\Hp^n)$. Moreover the constant $(n-2)^{-\frac{p+2}{p}} \overline{S}_{n,p}$ is sharp.
\la{thm:abcde}
\end{theorem}
The function $Y(t)$ above can be compared with the logarithmic function $X(t)$ near zero; see Theorem \ref{thm:p_h_s} for a precise statement.

In case $n=3$ we actually have sharpness of the constant with a logarithmic factor:
\begin{theorem}
{\bf (Poincar\'{e}-Hardy-Sobolev inequality)}
There exists an $\overline{\alpha}_3 >0$ such that  for all $0<\alpha \leq\overline{\alpha}_3$ and all $v\in \cic(\Hp^3)$ there holds
\be
\int_{\Hp^3}|\nabla_{\Hp^3} v|^2dV \geq  \int_{\Hp^3}v^2 dV +
\frac{1}{4} \int_{\Hp^3}\frac{v^2}{\sinh^2\rho}dV  
+ S_{3,p} \Big( \int_{\Hp^n} (\sinh\rho)^{\frac{p-6}{2}} 
X^{\frac{p+2}{2}}\big( \alpha\tanh(\rho/2)\big)|v|^{p} dV \Big)^{2/p}
\la{ell}
\ee
Moreover the constant  $S_{3,p}$ is sharp for all $0<\alpha\leq \overline{\alpha}_3$.
\la{thm:abcde3}
\end{theorem}

We next consider analogous inequalities in the case where the singularity is placed on the boundary of a bounded Euclidean domain $\Omega$ satisfying
an exterior ball condition.
Such Hardy-Sobolev inequalities have recently been obtained in \cite{BFT}.
Our aim here is to provide estimates for the Sobolev constant.
For $n\geq 3$ and $ 0\leq \gamma <n/2$ we denote by $S^*_{n,\gamma}$ the best constant for the inequality
\[
\Int_{\R^n_+}|\nabla u|^2dx - \gamma(n-\gamma)\Int_{\R^n_+}\frac{u^2}{|x|^2}dx 
\geq S^{*}_{n,\gamma}  
\Big( \int_{\R^n_+} |u|^{2^*}dx  \Big)^{2/2^*}   , \quad u\in \cic(\R^n_+).
\]
We then have
\begin{theorem}
\label{c}
{\bf (Hardy-Sobolev inequality)} 
Let $\Omega\subset\R^n$, $n\geq 3$, be a bounded domain with $0 \in \partial \xO$  and let $D=\sup_{\Omega}|x|$.
Assume that $\Omega$ satisfies an exterior ball condition
at zero, that is there exists a ball $B_{\rho} \subset \cC \overline{\xO} $. 
Then given $\gamma\in [0,n/2)$ there exist $r_{n,\gamma}$ and $\alpha_{n,\gamma}^*$  in $(0,1)$ both depending only on $n$ and $\gamma$ such that,
if the radius $\rho$  of the exterior ball satisfies  $\rho\geq D/ r_{n,\gamma}$ then for all $0<\alpha\leq\alpha_{n,\gamma}^*$ there holds
\[
\int_{\xO} |\nabla u|^2 dx \geq   \frac{n^2}{4}  \int_{\xO} \frac{u^2}{|x|^2} dx \\
+  (n-2\gamma)^{- \frac{2(n-1)}{n}}  S^*_{n,\gamma} \bigg( \int_{\xO}  X^{\frac{2n-2}{n-2}}|u|^{\frac{2n}{n-2}}     dx \bigg)^{\frac{n-2}{n}},
\]
for all $u \in C^{\infty}_{c}(\xO)$; here $X=X(\alpha |x|/ D)$.
\end{theorem}
In Theorem \ref{c_p} we establish a more general result where the RHS involves a weighted $L^p$ norm, $2<p\leq 2^*$.

The structure of the article is simple: in Section \ref{section:2} we study Hardy-Sobolev inequalities in Euclidean space with an interior singularity, in Section \ref{section:3} we study Hardy-Sobolev inequalities in hyperbolic space, whereas in Section \ref{section_bps} we study  Hardy-Sobolev inequalities when the singularity is placed on the boundary of a Euclidean domain $\Omega$.

\section{Hardy-Sobolev inequalities on Euclidean space}
\la{section:2}

In this section we establish improved Hardy-Sobolev inequalities in the Euclidean space with an interior point singularity. Our first result reads
\begin{theorem}
Let $n\geq 3$,  $2<p\leq 2^*$ and $0\leq \theta <1/2$.
For any $0<\alpha\leq\alpha_n$ and for any $u\in \cic(B_1)$ there holds
\bea
&&\hspace{-2cm}\int_{B_1}|\nabla u|^2dx -\Big( \frac{n-2}{2}\Big)^2 \int_{B_1}\frac{u^2}{|x|^2}dx
-\theta(1-\theta) \int_{B_1}\frac{u^2}{|x|^2} X(\alpha |x|)^2 dx \nonumber\\[0.2cm]
&\geq&
\Big( \frac{1-2\theta}{n-2}\Big)^{\frac{p+2}{p}}S_{n,p}
 \Big( \int_{B_1} |x|^{\frac{p(n-2)}{2} -n}  X(\alpha |x|)^{\frac{p+2}{2}} |u|^{p}  dx \Big)^{2/p} .
\la{mel1}
\eea
Moreover the constant $( \frac{1-2\theta}{n-2})^{\frac{p+2}{p}}S_{n,p}$ is sharp for any choice of the parameters.
\la{thm:2019}
\end{theorem}
{\em Proof.} Let $\tau_{n,p,\theta,\alpha}$ denote the best Sobolev constant for the above inequality. Setting
\[
u(x) =|x|^{-\frac{n-2}{2}} X(\alpha |x|)^{-\theta}v(x)
\]
and using polar coordinates we find
\[
\tau_{n,p,\theta,\alpha} =\inf \frac{\Int_0^1 \Int_{\sph} r X(\alpha r)^{-2\theta}\Big(  v_r^2 +\frac{1}{r^2}|\nabla_{\omega}v|^2\Big) dS\, dr }
{ \bigg(  \Int_0^1 \Int_{\sph} r^{-1}X(\alpha r)^{\frac{p+2}{2}-\theta p}|v|^{p}dS \, dr \bigg)^{2/p}} .
\]
Setting
\[
t= X(\alpha r)^{2\theta-1}  \; , \qquad v(r,\omega) =w(t,\omega),
\]
we then obtain
\be
\la{cp}
\tau_{n,p,\theta, \alpha} =(1-2\theta)^{\frac{p+2}{p}}\inf \frac{\Int_{X(\alpha)^{2\theta -1}}^{\infty} \Int_{\sph} 
\Big(  w_t^2 + (1-2\theta)^{-2} t^{\frac{4\theta}{1-2\theta}} |\nabla_{\omega}w|^2\Big) dS\, dt }
{ \bigg(  \Int_{X(\alpha)^{2\theta -1}}^{\infty} \Int_{\sph}  t^{-\frac{p+2}{2}} |w|^{p}dS \, dt \bigg)^{2/p}}
\ee
On the other hand we have for any $R>0$,
\bean
S_{n,p} &=&  \inf \frac{ \Int_{\R^n}|\nabla u|^2 dx}
{ \Big( \Int_{\R^n}  |x|^{\frac{p(n -2)}{2} -n }|u|^{p} dx \Big)^{2/p}  }   \\
&=&  \inf \frac{ \Int_{B_R}|\nabla u|^2 dx}
{ \Big( \Int_{B_R}  |x|^{\frac{p(n -2)}{2} -n }|u|^{p} dx \Big)^{2/p}  }  \\
&=& \inf \frac{ \Int_{0}^{R}\Int_{\sph} r^{n-1} \Big(  u_r^2 +\frac{1}{r^2}|\nabla_{\omega}u|^2 
\Big) dS\, dr}{ \Big( \Int_{0}^{R}\Int_{\sph} r^{\frac{p(n-2)}{2}-1 }|u|^{p}  dS\, dr \Big)^{2/p}  }.
\eean
Making the change of variables
\[
t=r^{-(n-2)} \; , \qquad u(r,\omega)=v(t,\omega)
\]
we easily arrive at
\be
(n-2)^{- \frac{p+2}{p}}  S_{n,p} =  \inf \frac{ \Int_{R^{2-n}}^{\infty}\Int_{\sph} \Big(  v_t^2 +\frac{1}{(n-2)^2t^2}|\nabla_{\omega}v|^2 
\Big) dS\, dt}{ \Big( \Int_{R^{2-n}}^{\infty}\Int_{\sph} t^{-\frac{p+2}{2}}|v|^{p}  dS\, dt \Big)^{2/p}  }.
\la{rain111}
\ee
We now choose $R$ so that $R^{2-n}= X(\alpha)^{2\theta -1}$. We also note that for all $0<\alpha\leq \alpha_n$ and $0\leq \theta<1/2$ there holds
$ X(\alpha)^{-1} \geq  (1-2\theta)/(n-2)$. Therefore
\[
(1-2\theta)^{-2} t^{\frac{4\theta}{1-2\theta}}  \geq \frac{1}{(n-2)^2t^2} \; , \qquad \mbox{ for all }  t\geq X(\alpha)^{2\theta -1} .
\]
Hence (\ref{mel1}) follows from (\ref{cp}) and (\ref{rain111}).
The sharpness follows from the fact that for any $R>0$,
\[
S_{n,p}=  \inf_{{\substack{u\in C^{\infty}_c(B_R) \\   \mbox{ \scriptsize{$u$ radial}} }}}  \frac{ \Int_{B_R}|\nabla u|^2 dx}
{ \Big( \Int_{B_R}  |x|^{\frac{p(n -2)}{2} -n }|u|^{p} dx \Big)^{2/p}  } \, .
\]
$\hfill\Box$


\

Let $\theta\in (0,2)$ and $R > 1$. For $r\in (0,1)$ we define the function
\bea
B(r) &=& \frac{1}{ (R^{\theta} -r^{\theta})^2 \Big(  1+ \int_r^1 \frac{ds}{s(R^{\theta} -s^{\theta})^2} \Big)} \nonumber \\
&=&
\frac{\theta R^{2\theta}} { \Big(  \theta R^{2\theta}  -\ln(R^{\theta} -1) + \frac{R^{\theta}}{R^{\theta} - 1}  -\ln( \frac{r^{\theta}}{ R^{\theta}-r^{\theta}})
-  \frac{R^{\theta}}{R^{\theta} -r^{\theta}} \Big)  (R^{\theta} -r^{\theta})^2} .
\label{beta}
\eea
\begin{lemma}
Let $\theta\in (0,2)$ and $R > 1$. Let $\alpha$ and $\beta$ be defined by
\bean
&& -\ln\alpha = R^{2\theta} -1 + \int_0^1 \frac{s^{\theta-1}(  2R^{\theta} -s^{\theta} )}{ (R^{\theta} -s^{\theta})^2}ds \; , \\
&& -\ln\beta =(R^{\theta} -1)^2 -1  >-1 \, .
\eean
Then there holds
\be
X(\alpha r) \leq B(r) \leq X(\beta r)  \; , \qquad r\in (0,1).
\la{sig}
\ee
\la{lem:bx}
\end{lemma}
{\em Proof.}  We first note that (\ref{sig}) is written equivalently as
\[
 -\ln\beta \leq  \frac{1}{B(r)} - \frac{1}{X(r)} \leq -\ln\alpha  \; , \qquad r\in (0,1),
\]
We thus consider the function
\[
g(r) = \frac{1}{B(r)} -\frac{1}{X(r)} \; , \qquad r\in (0,1),
\]
We have
\[
B(r) = \frac{1}{ t(R^{\theta} -r^{\theta})^2},
\]
therefore
\[
g(r) = ( R^{\theta} -r^{\theta})^2\Big(  1+ \int_r^1 \frac{ds}{s(R^{\theta} -s^{\theta})^2} \Big) -1+ \ln r .
\]
From this easily follows that $g$ is decreasing. Now, we can write $g(r)$ in an equivalent way as
\[
g(r) =-1 + (R^{\theta} -r^{\theta})^2 \Big(   1+ \frac{1}{R^{2\theta}}
 \int_r^1 \frac{s^{\theta-1} (2R^{\theta} -s^{\theta})}{ (R^{\theta}-s^{\theta})^2} ds  \Big)
\]
and therefore
\[
g(0+ ) = R^{2\theta} -1 +  \int_0^1 \frac{s^{\theta-1} (2R^{\theta} -s^{\theta})}{ (R^{\theta}-s^{\theta})^2}ds  \; .
\]
The result now follows from the monotonicity of $g$. $\hfill\Box$

We note that if we choose $\theta=0$ in Theorem \ref{thm:2019} we obtain inequality (\ref{intro:mel1}). For our purposes we shall also need
an improvement of (\ref{intro:mel1}) which we believe is of independent interest and reads as follows:
\begin{theorem}
Let $n\geq 3$, $2< p\leq 2^*$, $\theta\in (0,2)$.  We define $R>1$ and $\alpha_{n,\theta} <1$ by
\[
R^{\theta} = 1+ \frac{1}{\sqrt{n-2}}  \;\; , \qquad  -\ln\alpha_{n,\theta}= 
R^{2\theta} -1 + \int_0^1 \frac{s^{\theta-1}(  2R^{\theta} -s^{\theta} )}{ (R^{\theta} -s^{\theta})^2}ds
\]
Then for all $0<\alpha\leq \alpha_{n,\theta}$ and for all $u\in \cic(B_1)$ there holds
\begin{eqnarray}
&&\hspace{-3cm}\int_{B_1}|\nabla u|^2 dx  -\Big(  \frac{n-2}{2}\Big)^2 \int_{B_1}\frac{u^2}{|x|^2} dx -\theta^2 \int_{B_1} \frac{ u^2}{ |x|^{2-\theta} (R^{\theta} -|x|^{\theta})} dx \nonumber \\[0.2cm]
&\geq&  (n-2)^{-\frac{p+2}{p}}S_{n,p} \bigg(  \int_{B_1}  |x|^{\frac{p(n -2)}{2} -n} X^{\frac{p+2}{2}}(\alpha |x|) |u|^{p} dx  \bigg)^{2/p}.
\la{pin}
\end{eqnarray}
Moreover the constant $(n-2)^{-\frac{p+2}{p}}S_{n,p}$ is sharp.
\la{thm7}
\end{theorem}
{\em Proof.} A simple computation shows that the function
\[
\psi(x) =|x|^{ -\frac{n-2}{2}} (R^{\theta} -|x|^{\theta})
\]
satisfies
\[
- \frac{ \Delta\psi}{\psi} =\Big(  \frac{n-2}{2}\Big)^2 \frac{1}{|x|^2} + \frac{ \theta^2}{ |x|^{2-\theta} (R^{\theta} -|x|^{\theta})} .
\]
Making the change of variables $u=\psi v$ and using Lemma \ref{lem:bx} it then follows that in order to prove (\ref{pin}) it is enough to establish that
\be
\int_{B_1} |x|^{2-n} (R^{\theta} -|x|^{\theta})^2 |\nabla v|^2 dx \geq  
(n-2)^{-\frac{p+2}{p}}S_{n,p} \bigg(  \int_{B_1} |x|^{-n} (R^{\theta} -|x|^{\theta})^{p}
B(|x|)^{\frac{p+2}{2}} |v|^{p} dx \bigg)^{2/p}
\la{pin1}
\ee
for all $v\in \cic(B_1)$. 

Let $\tau_{n,\theta,p}$ denote the best constant for (\ref{pin1}), so that
\bean
\tau_{n,\theta,p}& =&\inf\frac{\Int_{B_1} |x|^{2-n} (R^{\theta} -|x|^{\theta})^2 |\nabla v|^2 dx}
{ \bigg(  \Int_{B_1} |x|^{-n} (R^{\theta} -|x|^{\theta})^{p}
B(|x|)^{\frac{p+2}{2}} |v|^{p} dx \bigg)^{2/p}} \\
&=& \inf \frac{ \Int_0^1 \Int_{\sph}  r (R^{\theta} -r^{\theta})^2 \Big(  v_r^2  + \frac{1}{r^2} |\nabla_{\omega}v|^2 \Big) dS \, dr}
{ \Big( \Int_0^1 \Int_{\sph} r^{-1} (R^{\theta} -r^{\theta})^{p}  B(r)^{\frac{p+2}{2}} |v|^{p} dS \, dr \Big)^{2/p}}.
\eean
We change variables setting
\be
v(r,\omega) =w(t,\omega) \; , \qquad t= 1+ \int_r^1 \frac{ds}{s(R^{\theta} -s^{\theta})^2}.
\la{kof}
\ee
It is easily seen that $r\mapsto t$ is a strictly decreasing map that maps $(0,1)$ onto $(1,+\infty)$ and
in addition we have
\[
\tau_{n,\theta,p} =\inf  \frac{  \Int_1^{\infty} \Int_{\us} \Big( w_t^2 + (R^{\theta}-r^{\theta})^4 |\nabla_{\omega}w|^2 \Big) dS \, dt }
{ \Big(  \Int_1^{\infty} \Int_{\us}   t^{-\frac{p+2}{2}}  |w|^p \, dS \, dt\Big)^{2/p}}
\]
We claim that
\be
(R^{\theta} -r^{\theta})^4 \geq \frac{1}{(n-2)^2t^2}  \; ,  \qquad r\in (0,1).
\la{xatz}
\ee
This is written equivalently as
\be
(R^{\theta} -r^{\theta})^2 (n-2)t \geq 1   \; ,  \qquad r\in (0,1),
\la{333}
\ee
where $t=t(r)$ is given by (\ref{kof}). Hence (\ref{xatz}) follows by noting that
the LHS of (\ref{333}) is a decreasing function of $r$, being the product of two positive decreasing functions. Therefore
\[
\tau_{n,\theta,p} \geq \inf  \frac{  \Int_1^{\infty} \Int_{\us} \Big( w_t^2 + \frac{1}{(n-2)^2t^2} |\nabla_{\omega}w|^2 \Big) dS \, dt }
{ \Big(  \Int_1^{\infty} \Int_{\us}   t^{-\frac{p+2}{2}}  |w|^p \, dS \, dt\Big)^{2/p}}
= (n-2)^{-\frac{p+2}{p}}S_{n,p} ,
\]
by (\ref{rain111}). 

The sharpness  of the constant of Theorem \ref{thm:2019} for the choice $\theta=0$ implies the sharpness of the constant
$(n-2)^{-\frac{p+2}{p}}S_{n,p}$ in (\ref{pin}). $\hfill\Box$


\section{Hardy-Sobolev inequalities on hyperbolic space}
\la{section:3}

In this section we study Hardy-Sobolev inequalities on the hyperbolic space $\Hp^n$. There are two standard models for $\Hp^n$. The first one is the unit ball model, where the unit ball $B_1$ is equipped with the Riemannian metric
\[
ds^2 = \Big(  \frac{1-|x|^2}{2} \Big)^{-2}dx^2 .
\]
Under this model we have
\[
 |\nabla_{\Hp^n}v|^2  =\Big(  \frac{1-|x|^2}{2} \Big)^2  |\nabla_{\R^n} v|^2 \quad , \quad \qquad 
  dV =\Big(  \frac{1-|x|^2}{2} \Big)^{-n} dx \, .
\]
Denoting by $\rho(x)$ the distance of $x\in B_1$ to the origin we then have
\[
 \rho(x) =\ln\Big( \frac{1+|x|}{1-|x|} \Big)  .
\]

We shall also use the half-space model of $\Hp^n$, namely $\R^n_+$ equipped with the Riemannian metric
\[
ds^2 =\frac{dx^2}{x_n^2}.
\]
Under this model we have
\[
 |\nabla_{\Hp^n}v|^2  = x_n^2 |\nabla_{\R^n} v|^2  \;\;\; , \qquad
  dV = x_n^{-n} dx \; .
\]

\subsection{Hardy-Sobolev inequalities}

In this subsection we prove Hardy-Sobolev inequalities that are analogues to those of the Euclidean case for an interior point singularity.
We start with
\begin{theorem} {\bf (Hardy-Sobolev inequality I)}
Let $n\geq 3$ and $2<p\leq 2^*$. For any $0<\alpha \leq\alpha_n$ there holds
\be
\int_{\Hp^n}|\nabla_{\Hp^n} v|^2dV -\Big( \frac{n-2}{2}\Big)^2 \int_{\Hp^n}\frac{v^2}{\rho^2}dV 
\geq
(n-2)^{-\frac{p+2}{p}}S_{n,p}\Big( \int_{\Hp^n}  (\sinh \rho)^{\frac{p(n-2)}{2}-n} X^{\frac{p+2}{2}} |v|^{p} dV \Big)^{2/p} ,
\la{eq:10}
\ee
for all $v\in\cic(\Hp^n)$; here $X=X(\alpha \tanh (\rho/2))$. Moreover the constant $(n-2)^{-\frac{p+2}{p}}S_{n,p}$ is sharp for all values of $\alpha\in (0,\alpha_n]$.
\la{thm1}
\end{theorem}
{\em Proof.} We use the ball model $B_1$ for $\Hp^n$, taking the centre of $B_1$ to correspond to the point $x_0$ where distance is taken from.
Using this model the required inequality (\ref{eq:10}) takes the form
\bea
&&\hspace{-1.5cm}\int_{B_1} \Big(  \frac{1-|x|^2}{2} \Big)^{2-n} |\nabla v|^2 dx   -\Big( \frac{n-2}{2}\Big)^2 \int_{B_1} \Big(  \frac{1-|x|^2}{2} \Big)^{-n} 
\frac{v^2}{\ln^2\Big( \frac{1+|x|}{1-|x|} \Big)} dx  \nonumber \\
&\geq& (n-2)^{-\frac{p+2}{p}}S_{n,p}
\bigg( \int_{B_1} |x|^{ \frac{p(n-2)}{2} -n}  \Big(  \frac{1-|x|^2}{2} \Big)^{-\frac{p(n-2)}{2}}  X^{\frac{p+2}{2}}(\alpha |x|)
|v|^{p} dx \bigg)^{2/p} ,
\la{b1}
\eea
for all $v\in\cic(B_1)$.

To prove (\ref{b1}) we use Theorem \ref{thm:2019} with the choice $\theta=0$. Making the change of variables
\[
u(x) =\Big(  \frac{1-|x|^2}{2} \Big)^{-\frac{n-2}{2}}  v(x).
\]
in (\ref{mel1}) we obtain
\bean
&&\hspace{-1cm}\int_{B_1} \Big(  \frac{1-|x|^2}{2} \Big)^{2-n} |\nabla v|^2 dx  -\frac{n(n-2)}{4} 
\int_{B_1} \Big(  \frac{1-|x|^2}{2} \Big)^{-n}   v^2 dx 
 -\Big( \frac{n-2}{2}\Big)^2 \int_{B_1} \Big(  \frac{1-|x|^2}{2} \Big)^{2-n}  \frac{v^2}{|x|^2} dx \nonumber  \\
&\geq & (n-2)^{-\frac{p+2}{p}}S_{n,p}
\bigg( \int_{B_1} |x|^{ \frac{p(n-2)}{2} -n}  \Big(  \frac{1-|x|^2}{2} \Big)^{-\frac{p(n-2)}{2}}  X^{\frac{p+2}{2}}(\alpha |x|)
|v|^{p} dx \bigg)^{2/p}
\eean
Therefore to prove (\ref{b1}) it is enough to establish
\bean
&&\frac{n(n-2)}{4} 
 \Big(  \frac{1-|x|^2}{2} \Big)^{-n}    
 + \Big( \frac{n-2}{2}\Big)^2  \Big(  \frac{1-|x|^2}{2} \Big)^{2-n}  \frac{1}{|x|^2}  \\
&\geq&  
\Big( \frac{n-2}{2}\Big)^2   \Big(  \frac{1-|x|^2}{2} \Big)^{-n} \frac{1}{  \big( \ln \big[\frac{1+|x|}{1-|x|} \big] \big)^2},
\eean
for all $x\in B_1$.
This is written equivalently as
\[
\Big( \ln \Big[\frac{1+|x|}{1-|x|}\Big] \Big)^2 \bigg\{   \Big(  \frac{1-|x|^2}{2} \Big)^{2} +\frac{n}{n-2}|x|^4 
+ \frac{2n}{n-2}  \Big(  \frac{1-|x|^2}{2} \Big) |x|^2    \bigg\}  \geq |x|^2 .
\]
This follows setting $t=|x|$ and using the elementary inequality
\[
\ln \frac{1+t}{1-t} \geq 2t , \qquad 0<t<1.
\]
The sharpness of the constant is a consequence of the sharpness of the constant of Theorem \ref{thm:2019}.
$\hfill\Box$


We next have 

{\bf\em Proof of Theorem \ref{thm:intro1}.} We first make the substitution
\[
u(x) =(\sinh \rho)^{-\frac{n-2}{2}} v(x).
\]
We then have
\bean
\int_{\Hp^n}|\nabla_{\Hp^n} u|^2dV &=& \int_{\Hp^n} (\sinh \rho)^{-n+2} |\nabla_{\Hp^n} v|^2dV +
\Big( \frac{n-2}{2}\Big)^2 \int_{\Hp^n} (\sinh \rho)^{-n}  v^2 dV \\[0.2cm]
&& +\frac{n(n-2)}{4} \int_{\Hp^n} (\sinh \rho)^{-n+2}  v^2 dV \, ,
\eean
hence the required inequality (\ref{101}) becomes
\[
\int_{\Hp^n} (\sinh \rho)^{-n+2}  |\nabla_{\Hp^n} v|^2dV \geq
(n-2)^{-\frac{p+2}{p}}S_{n,p}
\Big( \int_{\Hp^n}  (\sinh \rho)^{-n}   X^{\frac{p+2}{2}} \big(\alpha \tanh (\rho/2)\big)
|v|^{p} dV \Big)^{2/p} .
\]
Using spherical hyperbolic coordinates around $x_0$ this is written
\bea
&& \hspace{-2.5cm} \int_0^{\infty}\int_{\sph} (\sinh \rho) \Big( v_{\rho}^2 +\frac{1}{\sinh^2\rho}  |\nabla_{\omega}v|^2 \Big) dS \, d\rho 
\la{pier} \\
&\geq&
(n-2)^{-\frac{p+2}{p}}S_{n,p}
\Big( \int_0^{\infty} \int_{\sph}  (\sinh \rho)^{-1}  X^{\frac{p+2}{2}} \big(\alpha \tanh (\rho/2)\big)
|v|^{p} dS \, d\rho \Big)^{2/p}. \nonumber
\eea
To prove (\ref{pier}) we change variables setting
\[
\frac{1}{t} =X\big(\alpha \tanh (\rho/2) \big) \; , \qquad 
v(\rho,\omega) =w(t,\omega) .
\]
We then have
\bean
\int_0^{\infty}\int_{\sph} (\sinh \rho) \Big( v_{\rho}^2 +\frac{1}{\sinh^2\rho}  |\nabla_{\omega}v|^2 \Big) dS \, d\rho  &=&
\int_{X^{-1}(\alpha)}^{\infty}\int_{\sph} \Big( w_t^2 +|\nabla_{\omega}w|^2 \Big) dS \, dt \\
&\geq& \int_{X^{-1}(\alpha)}^{\infty}\int_{\sph} \Big( w_t^2 +\frac{1}{(n-2)^2t^2}|\nabla_{\omega}w|^2 \Big) dS \, dt \, ,
\eean
since $X^{-1}(\alpha) \geq 1/(n-2)$. We also have
\[
\int_0^{\infty} \int_{\sph}  (\sinh \rho)^{-1}   X^{\frac{p+2}{2}}\big(\alpha \tanh (\rho/2)\big) 
|v|^{p} dS \, d\rho =  \int_{X^{-1}(\alpha)}^{\infty}\int_{\sph}  t^{-\frac{p+2}{2}}|w|^{p} dS \, dt \,.
\]
Inequality (\ref{pier}) now follows from the last two relations together with (\ref{rain111}). The sharpness of the constant 
$(n-2)^{-\frac{p+2}{p}}S_{n,p}$ is a consequence of the sharpness of the constant of Theorem \ref{thm1}.
$\hfill\Box$

\

{\bf Remark.}
Since
\[
\Big( \frac{n-2}{2}\Big)^2 \frac{1}{\sinh ^2\rho} +\frac{n(n-2)}{4} \geq  \Big( \frac{n-2}{2}\Big)^2 \frac{1}{\rho^2}  \;\; , \qquad \rho>0,
\]
one can obtain inequality (\ref{eq:10}) as a consequence of inequality (\ref{101}). However the sharpness of the constant in (\ref{eq:10}) does not follow.




\subsection{Poincar\'{e}-Hardy-Sobolev inequalities}

In this subsection we give the proof of Theorem \ref{thm:abcde} and other related results.

We recall the definition of the hypergeometric function
\[
F(a,b,c,z )= \frac{\Gamma(c)}{\Gamma(a)\Gamma(b)} \sum_{n=0}^{\infty}\frac{\Gamma(a+n)\Gamma(b+n)}{\Gamma(c+n)} \frac{z^n}{n!}  \; , \qquad |z|<1. 
\]
We refer to \cite[Section15]{AS} for various properties of the hypergeometric functions.

We shall use three specific hypergeometric functions. The first one is
\[
F(z) =F\big( \frac{1}{2} , \frac{1}{2} , 1 , z\big)  \; , \quad |z|<1 \, .
\]
We then have
\begin{lemma}
The problem
\be
\darr{g'' + \Frac{1}{4\sinh^2 t}g=0 \; ,}{ \qquad t>0}
{\lim_{t\to +\infty}g(t)=1 \, .}{}
\la{g}
\ee
has a unique  solution $g(t)$. Moreover the solution is  positive, strictly increasing and is given by
\be
g(t)=
\left\{
\begin{array}{ll}
{ (e^{2t}-1)^{1/2}F( 1-e^{2t})  \Big( { 1 +\Frac{1}{\pi}}
 \Int_{-1}^{ 1-e^{2t} } \frac{ds}{s(s-1)F(s)^2} \Big) ,} & {\;\;  0<t< \ln\sqrt{2} \,  ,}  \\[0.3cm]
{F  \big( \frac{1}{1-e^{2t}}  \big) \; ,}  & { \;\;  t\geq \ln\sqrt{2} \, .}
\end{array} \right.
\la{g:formula}
\ee
\la{lem:g}
\end{lemma}
{\em Proof.} We shall first prove that if a solution of (\ref{g}) exists then it
is positive in $(0,+\infty)$. Suppose to the contrary that there exists a $\rho>0$  such that $g(\rho)=0$.
Then
\[
\frac{1}{4} \int_{\rho}^{\infty}  \frac{g^2}{\sinh^2 t}dt = - \int_{\rho}^{\infty}  g'' g \, dt  =  \int_{\rho}^{\infty}  (g')^2  dt  
\]
Since $t< \sinh t$ for $t>0$, this implies that $g$ is a minimizer for the standard Hardy inequality in $(\rho,+\infty)$, which is a contradiction since there is no minimizer. The monotonicity of $g(t)$ then also follows immediately.

We now change variables in (\ref{g}) setting
\[
 \xi = \frac{1}{1-e^{2t}} \;\; , \qquad \quad g(t) =A(\xi).
\]
Simple computations show that $g(t)$ is a solution of (\ref{g}) if and only if $A(\xi)$ is a solution of
\be
\xi(\xi-1)A'' +(2\xi-1)A' +\frac{1}{4} A =0  \; , \qquad \xi <0.
\la{110}
\ee
For $\xi \in (-1,0)$ one solution of (\ref{110}) is the hypergeometric function $F(\xi)$ which  is positive for $\xi\in (-1,0)$.
A second solution of (\ref{110}) for $\xi\in (-1,0)$ can be found by standard arguments; after simple computations we find that a second solution is
 \[
   F(\xi)  \int_{-1}^{-\xi} \frac{ds}{s(s-1)F(s)^2}  \;  , \qquad    -1< \xi<0 .
 \]
We note in particular that
\[
F(\xi)  \int_{-1}^{-\xi} \frac{ds}{s(s-1)F(s)^2}  \sim \ln(-\xi)  \; , \qquad \mbox{ as }\xi \to 0-  .
\]
Therefore the solution of (\ref{110}) in $(-1,0)$ with $A(0)=1$ is the function $F(\xi)$ and hence the solution of (\ref{g}) with $\lim_{t\to +\infty}g(t)=1$ is unique and is given by
\[
g(t)   =  F  \big( \frac{1}{1-e^{2t}}  \big) \; \;   , \qquad   t\geq \ln\sqrt{2}. 
\]
We want to describe the solution $g(t)$ for $t\in (0,\ln\sqrt{2})$ and to do so we extend the corresponding function $A(\xi)$ to the interval
$(-\infty,-1)$. For this we first note that the function
\[
B(\xi) = (-\xi)^{-1/2} A(\frac{1}{\xi})  \; , \qquad \xi<0,
\]
is a solution of (\ref{110}) if and only if $A(\xi)$ is a solution of (\ref{110}).  It follows that the extension of $A(\xi)$ has the form
\[
A(\xi) =  (-\xi)^{-1/2}F(\frac{1}{\xi}) \Big( c_1^*  + c_2^* 
 \Int_{-1}^{\frac{1}{\xi}   } \frac{ds}{s(1-s)F(s)^2} \Big)  \;\;  , \qquad \xi<-1 \, ,
\]
for suitable constants $c_1^*$, $c_2^*$.  The continuity of $A(\xi)$ at $\xi=-1$ gives $c_1^*=1$. Moreover, the differentiability of $A(\xi)$
at $\xi=-1$ gives
\[
c_2^* =   F(-1) \big(   4F'(-1) -F(-1) \big)  =-\frac{1}{\pi}.
\]
This concludes the proof. $\hfill\Box$

\

{\bf Remark.} We note for future use that the above proof implies in particular that the hypergeometric function $F(\xi)$ is positive and increasing in the interval $(-1,0)$.


In our next lemma we make use of the hypergeometric functions
\[
F_1(z) =F( \frac{n-1}{2}, \frac{n-1}{2} , n-1 , z) \; , \qquad  |z|<1 \, ,
\]
and
\[
F_2(z) =F(\frac{n-1}{2} , -\frac{n-3}{2} , 1 , z)   \; , \qquad |z|<1 \, .
\]

\begin{lemma}
Let $n\geq 4$. The problem
\be
\darr{h'' - \Frac{(n-1)(n-3)}{4\sinh^2  t}h=0}{ \; , \qquad  t>0,}{\lim_{t\to +\infty}h( t)=1 \, }{}
\la{h}
\ee
has a unique solution $h( t)$. Moreover the solution is positive, strictly decreasing and there exist constants $c_1^{\#}$, $c_2^{\#}$ with 
 $c_2^{\#} \neq 0$ so that
\be
h( t)=
\left\{
\begin{array}{ll}
{ (e^{2 t}-1)^{\frac{n-1}{2}}F_1( 1-e^{2 t}) \Big( c_1^{\#}  + c_2^{\#} 
 \Int_{-1}^{ 1-e^{2 t} } \frac{ds}{s^{n-1}(s-1)F_1(s)^2} \Big) ,} & {\;\;  0< t< \ln\sqrt{2} \,  ,}  \\[0.3cm]
{F_2 \big( \frac{1}{1-e^{2 t}}  \big) \; ,}  & { \;\;   t \geq \ln\sqrt{2} \, ,}
\end{array} \right.  \; .
\la{h:formula}
\ee
\la{lem:h}
\end{lemma}
{\em Proof.} The monotonicity and positivity of any solution $h(t)$ of (\ref{h})  follow from the differential equation by a simple argument.

We change variables setting
\[
 \xi = \frac{1}{1- e^{2 t}} \; \; , \qquad   h( t)=A(\xi).
\]
Equation (\ref{h}) then becomes
\be
\xi(\xi-1)A'' +(2\xi-1)A' -\frac{(n-1)(n-3)}{4} A =0  \; , \qquad \xi <0.
\la{210}
\ee
One solution of (\ref{210}) for $\xi \in (-1,0)$ is the hypergeometric function
$F_2(\xi)$ defined above.
A second solution of the ODE for $\xi\in (-1,0)$ is again found by standard arguments to be the function
 \[
   F_2(\xi)  \int_{-1}^{-\xi} \frac{ds}{s(s-1)F_2(s)^2}    \;  , \qquad    -1< \xi<0 .
 \]
which behaves like $\ln(-\xi)$ as $\xi\to 0-$.
Hence the solution of (\ref{210}) in $(-1,0)$ with $A(0)=1$ is the function $F_2(\xi)$ and therefore problem (\ref{h})
has a unique solution $h(t)$ which for $t\geq \ln\sqrt{2}$ is given by
\be
h(t)   =  F_2  \big( \frac{1}{1-e^{2t}}  \big) \; \;   , \qquad   t\geq \ln\sqrt{2}. 
\la{h:f}
\ee
We next describe the solution $h(t)$ for $t\in (0,\ln\sqrt{2})$ and to do so we extend the corresponding function $A(\xi)$ to the interval
$(-\infty,-1)$. For this we first note that the transformation
\[
A(\xi) =(-\xi)^{-\frac{n-1}{2}} B(\frac{1}{\xi}) \; , \qquad \xi<0 \, ,
\]
transforms a solution $A(\xi)$, $\xi <0$, of (\ref{210}) to a solution $B(w)$, $w<0$, of
\be
w(w-1)B'' +(n w +1-n)B'  +\Big(\frac{n-1}{2}\Big)^2 B =0    \; , \qquad w <0.
\la{211}
\ee
One solution of (\ref{211}) for $w\in (-1,0)$ is the hypergeometric function $F_1(w)$ defined above. A second solution is the function
\[
F_1(w)\Int_{-1}^{w } \frac{ds}{s^{n-1}(s-1)F_1(s)^2}.
\]
We thus conclude that there exist  $c_1^{\#}$, $c_2^{\#}$ so that
\[
A(\xi)= (-\xi)^{-\frac{n-1}{2}}F_1(\frac{1}{\xi}) \Big( c_1^{\#}  + c_2^{\#} 
 \Int_{-1}^{\frac{1}{\xi}   } \frac{ds}{s^{n-1}(s-1)F_1(s)^2} \Big) \;\; , \quad \xi<-1 \, .
\]
By the differentiability of $A(\xi)$ at $\xi=-1$ we obtain
\[
c_1^{\#} =\frac{F_2(-1)}{F_1(-1)} \;\; , \qquad  c_2^{\#} =2(-1)^{n-1}F_1(-1)F_2(-1)\Big(  \frac{F_1'(-1) }{F_1(-1)} + \frac{F_2'(-1) }{F_2(-1)}
-\frac{n-1}{2} \Big).
\]
It remains to prove that $c_2^{\#}\neq 0$. By (\ref{h:f}) and the fact that $h(t)$ is strictly decreasing we have $F_2'(-1)/F_2(-1)<0$ and hence it suffices to establish that
\be
\frac{F_1'(-1) }{F_1(-1)} \leq \frac{n-1}{4}.
\la{com}
\ee
Let as define the function
\[
q(w) = \frac{F_1'(w) }{F_1(w)}  \; , \quad  -1\leq w \leq 0.
\]
This function is a solution of
\[
q' =\frac{1}{w(1-w)} \bigg[  \Big(  \frac{n-1}{2}\Big)^2  -( n-1-nw )q -w(1-w)q^2  \bigg]  \; , \qquad -1\leq w<0.
\]
Moreover an elementary computation gives
\[
q'(0) =\frac{(n-1)(3n+1)}{16n} >0,
\]
which implies in particular that $q(w)<q(0)=(n-1)/4$ for $w<0$ close enough to zero.
Since
\[
0< \frac{1}{w(1-w)} \bigg[  \Big(  \frac{n-1}{2}\Big)^2  -\Big( n-1-n w  \Big)\frac{n-1}{4} -w(1-w)\Big( \frac{n-1}{4} \Big)^2  \bigg]
= \frac{(n-1)(3n+1+(n-1)w)}{16(1-w)}
\]
for $-1 \leq w<0$, a standard ODE comparison argument yields $q(w)<  (n-1)/4$, $-1\leq w <0$, and (\ref{com}) follows.  $\hfill\Box$


\begin{lemma}
The functions $g(t)$ and $h( t)$ satisfy the following asymptotic formulas as $t \to 0+ \, $:
\bea
&\ia &  g(t) =(2t)^{1/2}\big( -\frac{1}{\pi} \ln (2t) +B +O(t\ln t)  \, \big)    \nonumber \\
&& \mbox{where} \nonumber \\
&& \hspace{3cm}  B= 1+\frac{1}{\pi}\int_0^1 \frac{ 1- (1+t)F^2(-t)}{t(t+1)F^2(-t)}dt    \label{A} \\
&&   \mbox{and in particular }  1- \frac{1}{\pi} <B<1 . \nonumber  \\ [0.2cm]
& \ib & \mbox{If $n\geq 5$ then}  \nonumber  \\
&&  \hspace{2cm} h( t) =  \frac{c_2^{\#} (-1)^{n}}{n-2} (2 t)^{- \frac{n-3}{2} }  \Big(1 +\frac{(n-1)(n-3)}{24(n-4)} t^2 +O( t^3) \Big)  \nonumber  \\
&& \mbox{If $n=4$ then}  \nonumber \\
&&  \hspace{2cm} h( t) = \frac{c_2^{\#} }{2} (2 t)^{- \frac{1}{2} }  \Big(  1   -\frac{1}{8} t^2 \ln (2 t) +O( t^2) \Big) , \nonumber \\
&& \mbox{Here in both cases $c_2^{\#}$ is the non-zero coefficient of Lemma \ref{lem:h}.} \nonumber 
\eea
\la{lem:asympt_gh}
\end{lemma}
{\em Proof.} Part (i) follows from (\ref{g:formula}); we omit the details.
We next prove the double inequality for $B$.
To prove that $B<1$ it is enough to establish that $1- (1+t)F^2(-t) <0$ for $t\in (0,1)$ or equivalently $1 -(1-\xi)F(\xi)^2 <0$ for $\xi\in (-1,0)$. We have
\[
\frac{d}{d\xi} \big(1 -(1-\xi)F(\xi)^2\big) =F(\xi) \big( F(\xi) -2(1-\xi)F'(\xi) \big) =: F(\xi)Q(\xi).
\]
Since $F(0)=1$ the result will follow once we establish that $Q(\xi)>0$, $\xi\in(-1,0)$. Indeed, $Q(0)=1/2$ and
\[
Q'(\xi) =3F'(\xi) -2(1-\xi)F''(\xi) = \frac{1}{\xi}\big[ (2-\xi)F'(\xi)  +\frac{1}{4}F(\xi)   \big] <0 \; , \qquad -1<\xi<0.
\]
To prove that $B>1-1/\pi$ we first note that
\[
\frac{d}{d\xi}\big(  (1-\xi^2) F^2(\xi) \big) =2F(\xi) \big(   (1-\xi^2)F'(\xi) -2\xi F(\xi) \big) >0 \; , \qquad -1<\xi<0.
\]
Hence for $\xi\in (-1,0)$ we have $(1-\xi^2)F^2(\xi)<F^2(0)=1$; this implies
\[
 \frac{1 -(1-\xi)F^2(\xi)}{\xi(\xi-1)F^2(\xi)} >-1  \; , \qquad -1<\xi<0,
\]
and the result follows.

To prove (ii) we first recall (cf. (\ref{h:formula})) that
\[
h( t) =x^{\frac{n-1}{2}}F_1(-x)  \Big( c_1^{\#}  + c_2^{\#} 
 \Int_{-1}^{ -x } \frac{ds}{s^{n-1}(s-1)F_1(s)^2} \Big)  =:R(x) \; , \qquad\quad x=e^{2 t}-1 ,
\]
To estimate this we first note that
\[
F_1(x) = 1+\frac{n-1}{4}x + \frac{ (n-1)(n+1)^2}{32n}x^2 +O(x^3) \; , \qquad \mbox{ as }x\to 0,
\]
and therefore
\[
 \frac{1}{(s-1)F_1(s)^2} =-1 +\frac{n-3}{2}s -A_n s^2 +O(s^3),  \qquad \as s\to 0,
\]
where $A_n =(2n^3-15n^2+28n+1)/(16n)$. This implies that for small $x>0$ we have
\[
\int_{-1}^{-x} \frac{ds}{s^{n-1}(s-1)F_1(s)^2}=\darr{  \frac{(-1)^n}{n-2}x^{2-n} \Big(  1 +\frac{n-2}{2}x  +\frac{(n-2)A_n}{n-4}x^{2} +O(x^3)  \Big) ,}{ \mbox{ if }n\geq 5,}
{ \frac{1}{2}x^{-2}\Big(  1+x   -2A_4 x^2\ln x   +O(x^2)\Big),}{ \mbox{ if }n=4 \, .}
\]
Combining these we conclude that
\be
R(x)=\darr{ \frac{c_2^{\#}(-1)^n}{n-2}x^{-\frac{n-3}{2}}\Big( 1  + \frac{n-3}{4}x + B_nx^2 +O(x^3)  \Big),}
{\mbox{ if }n\geq 5,}{ \frac{c_2^{\#}}{2}x^{-\frac{1}{2}} \Big( 1 +\frac{1}{4}x -\frac{1}{32}x^2 \ln x  +O(x^2) \Big) ,  }{\mbox{ if }n=4 \, .}
\la{eq:q}
\ee
where
\[
B_n = \frac{(n-1)(n+1)^2}{32n} -\frac{(n-1)(n-2)}{8} +\frac{(n-2)A_n}{n-4} =\frac{ (n-3)(n-5)^2}{32(n-4)}.
\]
Setting $x=e^{2 t}-1$ we find
\[ 
\darr{1  + \frac{n-3}{4}x + B_nx^2 =1 + \frac{n-3}{2} t + \frac{ (n-3)^3}{8(n-4)} t^2 +O( t^3),}{ \mbox{ if }n\geq 5,}
{ 1 +\frac{1}{4}x -\frac{1}{32}x^2 \ln x = 1 +\frac{1}{2}  t -\frac{1}{8} t^2 \ln (2 t) +O( t^2),}{ \mbox{ if }n =4.}
\]
and
\[
x^{-\frac{n-3}{2}} =(2 t)^{-\frac{n-3}{2}} \Big( 1 -\frac{n-3}{2} t  + \frac{(n-3)(3n-11)}{24} t^2 +O( t^3) \Big).
\]
These together with (\ref{eq:q}) give the asymptotic formulas (ii).
$\hfill\Box$


\

{\bf\em Proof of Theorem \ref{thm:abcde}}. Using spherical coordinates in (\ref{intro:100}) we find that
\be
\overline{S}_{n,p} = \inf\frac{\Int_0^{\infty}\Int_{\sph}(\sinh\rho)^{n-1}\Big( u_{\rho}^2 + 
\frac{1}{ \sinh^2 \!\rho} |\nabla_{\omega} u|^2   -\Big( \frac{n-1}{2}\Big)^2 u^2   \Big)dS \, d\rho  }
{\Big(  \Int_0^{\infty}\Int_{\sph} (\sinh\rho)^{\frac{p(n-2)}{2} -1} |u|^{p} dS \, d\rho  \Big)^{2/p }}.
\la{eq:102}
\ee
Now let us define
\[
\phi(\rho) =(\sinh\rho)^{-\frac{n-1}{2}} h(\rho),
\]
where $h(\rho)$ is the function studied in Lemma \ref{lem:h} (in case $n=3$ we simply take $h(\rho)=1$). Simple computations then give that
$\phi(\rho)$ satisfies
\be
  \phi'' +(n-1)\coth \rho \,\, \phi +\Big( \frac{n-1}{2}\Big)^2 \phi =0  \; , \qquad \rho>0.
\la{eq:phi}
\ee
Setting $u=\phi \, w$ in (\ref{eq:102}) and using (\ref{eq:phi}) we then find that
\bea
\overline{S}_{n,p} & =& \inf\frac{\Int_0^{\infty}\Int_{\sph}  (\sinh\rho)^{n-1}\phi(\rho)^2  \Big( w_{\rho}^2 + 
\frac{1}{ \sinh^2 \!\rho} |\nabla_{\omega} w|^2   \Big)dS \, d\rho  }
{\Big(  \Int_0^{\infty}\Int_{\sph}   (\sinh\rho)^{\frac{p(n-2)}{2} -1} \phi(\rho)^{p}  |w|^{p} dS \, d\rho  \Big)^{2/p }}  \nonumber \\
& =& \inf\frac{\Int_0^{\infty}\Int_{\sph}  h(\rho)^2  \Big( w_{\rho}^2 + 
\frac{1}{ \sinh^2 \!\rho} |\nabla_{\omega} w|^2   \Big)dS \, d\rho  }
{\Big(  \Int_0^{\infty}\Int_{\sph}   (\sinh\rho)^{-  \frac{p+2}{2}    } h(\rho)^{p}  |w|^{p} dS \, d\rho  \Big)^{2/p }} .  \la{eq:103}  
\eea
Let us now define
\bea
\sigma_{n,p} =\inf\frac{\Int_{\Hp^n}|\nabla_{\Hp^n} u|^2dV - \Big( \frac{n-1}{2}\Big)^2 \int_{\Hp^n}u^2 dV -
\Big( \frac{n-2}{2}\Big)^2 \int_{\Hp^n}\frac{u^2}{\sinh^2\rho}dV}
{\Big( \Int_{\Hp^n} (\sinh\rho)^{\frac{p(n-2)}{2}-n} 
Y(\rho)^{\frac{p+2}{2}}|v|^{p} dV \Big)^{2/p}},
\la{eq:104}
\eea
the best constant for inequality (\ref{intro:102}). Using polar coordinates this is written as
\[
\sigma_{n,p} = \inf\frac{\Int_0^{\infty}\Int_{\sph}(\sinh t)^{n-1}\Big( u_t^2  +
\frac{1}{ \sinh^2 \! t} |\nabla_{\omega} u|^2   -\Big( \frac{n-1}{2}\Big)^2 u^2      -\Big( \frac{n-2}{2}\Big)^2 \frac{u^2}{  (\sinh t)^2}  \Big)dS \, dt  }
{\Big(  \Int_0^{\infty}\Int_{\sph} (\sinh \, t)^{\frac{p(n-2)}{2}-1} Y(t)^{\frac{p+2}{2}}  |u|^{p}   dS \, dt  \Big)^{2/p }}.
\]
We next define
\[
f(t) =(\sinh t)^{-\frac{n-1}{2}} g(t),
\]
where $g(t)$ is the solution of problem (\ref{intro:g}).
Simple computations then show that $f(t)$ is a positive solution to the equation
\be
f'' +(n-1)\coth  t \,  f'   +\Big(    \Big( \frac{n-1}{2}\Big)^2  
 + \Big( \frac{n-2}{2}\Big)^2 \frac{1}{\sinh^2 t} \Big)f =0  \; , \qquad t>0.
\la{eq:f}
\ee
Setting $u=fv$ in (\ref{eq:104}) and using (\ref{eq:f}) we obtain
\bea
\sigma_{n,p} &=& \inf\frac{\Int_0^{\infty}\Int_{\sph}(\sinh t)^{n-1}  f(t)^2\Big( v_t^2  +
\frac{1}{ \sinh^2 \! t} |\nabla_{\omega} v|^2  \Big)dS \, dt  }
{\Big(  \Int_0^{\infty}\Int_{\sph} (\sinh \, t)^{\frac{p(n-2)}{2}-1}  f(t)^{p} \, Y(t)^{\frac{p+2}{2}} |v|^{p}  dS \, dt  \Big)^{2/p }}  \nonumber \\
&=& \inf\frac{\Int_0^{\infty}\Int_{\sph}  g(t)^2\Big( v_t^2  +
\frac{1}{ \sinh^2 \! t} |\nabla_{\omega} v|^2  \Big)dS \, dt  }
{\Big(  \Int_0^{\infty}\Int_{\sph} (\sinh \, t)^{ -\frac{p+2}{2} }  g(t)^{p} \, Y(t)^{\frac{p+2}{2}} |v|^{p}  dS \, dt  \Big)^{2/p }}.
\la{eq:106}
\eea
We shall compare the expressions at the RHSs of (\ref{eq:103}) and (\ref{eq:106}) and in order to do so we change variables in (\ref{eq:106}) setting
\[
 \int_0^{\rho} \frac{dr}{h(r)^2} = \int_0^{t} \frac{ds}{g(s)^2}  \; , \qquad v(t,\omega) =w(\rho, \omega).
\]
(Here we note that both integrals are finite by the asymptotics of Lemma \ref{lem:asympt_gh}.) After some further computations we arrive at
\be
\sigma_{n,p}= (n-2)^{ -\frac{p+2}{p}} \inf
\frac{  \Int_{0}^{\infty}\Int_{\sph}h(\rho)^2 \Big(    w_{\rho}^2  + \Frac{g(t)^4}{ (\sinh t)^2 h(\rho)^4} |\nabla_{\omega}w|^2 \Big) dS \, d\rho}
{\Big(  \Int_0^{\infty}\Int_{\sph}   (\sinh\rho)^{-  \frac{p+2}{2}    } h(\rho)^{p}  |w|^{p} dS \, d\rho  \Big)^{2/p }}.
\la{eq:116}
\ee
Comparing (\ref{eq:103}) and (\ref{eq:116}) we conclude that in order to prove (\ref{intro:102}) it is enough to establish that 
\be
\frac{ g(t)^4}{  (\sinh t)^2  h(\rho)^4}    \geq 
\frac{ 1}{   (\sinh\rho)^{2}} .
\la{gh}
\ee
To see this we first recall from Lemmas \ref{lem:g} and \ref{lem:h} that $g(t)<1$ and $h(\rho)\geq 1$, so
\[
 \int_0^{t} \frac{ds}{g(s)^2} = \int_0^{\rho} \frac{dr}{h(r)^2}  <  \int_0^{\rho}  \frac{dr}{g(r)^2} \, , 
\]
and therefore $\rho>t$; this easily implies that
\be
\frac{t}{\rho} >  \frac{ \sinh  t}{ \sinh \rho} \, .
\la{par}
\ee
Therefore using the monotonicity of the functions $g(t)$ and $h(\rho)$,
\[
\frac{\rho}{h(\rho)^2} \geq   \int_0^{\rho} \frac{dr}{h(r)^2}  =  \int_0^{t} \frac{ds}{g(s)^2} > \frac{t}{g(t)^2} ,
\]
which, together with (\ref{par}), gives (\ref{gh}). This completes the proof of (\ref{intro:102}).

To prove the sharpness of the constant $\overline{S}_{n,p}$ we use decreasing rearrangements. It is easy to see (cf. \cite[Corollary 1]{Bae})  that the the infimum (\ref{eq:104}) remains the same if it is only taken amongst radial functions. This, together with  (\ref{eq:103}) and (\ref{eq:116}) implies the sharpness of the constant $\overline{S}_{n,p}$.
$\hfill\Box$

\begin{theorem}
{\bf (Poincar\'{e}-Hardy-Sobolev inequality II)}
Let $n\geq 3$ and $2<p\leq 2^*$.  There exists $0<\alpha_n <e$ which depends only on $n$ such that for
all $0<\alpha\leq \alpha_n$ and for all $ v\in \cic(\Hp^n)$ there holds
\bean
\int_{\Hp^n}|\nabla_{\Hp^n} v|^2dV &\geq& \Big( \frac{n-1}{2}\Big)^2 \int_{\Hp^n}v^2 dV +
\Big( \frac{n-2}{2}\Big)^2 \int_{\Hp^n}\frac{v^2}{\sinh^2\rho}dV   \\[0.2cm]
&& + (n-2)^{-\frac{p+2}{p}} \overline{S}_{n,p}\Big( \int_{\Hp^n} (\sinh\rho)^{\frac{p(n-2)}{2}-n} 
X^{\frac{p+2}{2}}\big(\alpha \tanh(\rho/2) \big)|v|^{p} dV \Big)^{2/p} .
\eean
\la{thm:p_h_s}
\end{theorem}
{\em Proof.} By Theorem \ref{thm:abcde} and by the monotonicity of $X$ it is enough to establish the existence of an $\alpha_n \in (0,e)$ such that
\be
Y(t)  \geq  X\big(  \alpha_n\tanh\frac{t}{2}  \big)  \; , \qquad t>0.
\la{yx}
\ee
By compactness, it is enough to prove that (\ref{yx}) is valid near zero and near infinity. 

{\em Case 1: Large $t>0$.} Let $\rho=\rho(t)$ be the function defined in (\ref{intro:fnr}).
We claim there exist $t_n>0$ and $c_n>0$ so that
\be
 \rho \leq t +c_n \; , \qquad \mbox{ for all }t\geq t_n \, .
\la{tn}
\ee
To prove this we first note that $F_2'(0)>0$ and hence there exists $A>0$ such that
\[
F_2(\xi)=F(\frac{n-1}{2} , -\frac{n-3}{2} , 1 , \xi) \leq 1-A\xi \; , 
\]
for all negative $\xi$ small enough in absolute value. It then follows from (\ref{h:formula}) that for large enough $\rho$ there holds
\[
\frac{1}{h(\rho)^2} = \frac{1}{ F_2^2 \big( \frac{1}{1-e^{2 \rho}}  \big)}
 \geq \frac{1}{ \big(  1+\frac{A}{e^{2\rho} -1}   \big)^2 }  =1 -\frac{A(2e^{2\rho}+A-2)}{(e^{2\rho}+A-1)^2}
\]
Hence if $\rho$ is large enough then
\be
 \int_0^{\rho} \frac{dr}{h(r)^2} =  \int_0^{\rho_0} \frac{dr}{h(r)^2} +  \int_{\rho_0}^{\rho} \frac{dr}{h(r)^2} \geq \rho -C_1 \, 
\la{h:lt}
\ee
for some $C_1>0$.

Analogous computations are valid for the function $g(t)$. We now use the estimate
\[
\frac{1}{g(t)^2} = \frac{1}{ F^2 \big( \frac{1}{1-e^{2 t}}  \big)} \leq  \frac{1}{ \big(  1+\frac{b}{1-e^{2t} }  \big)^2 }
= 1+ \frac{ b(e^{2t} -2-b)}{ (e^{2t} -1-b)^2},
\]
which is valid for some $b>0$ and large enough $t>0$. This leads to
\be
 \int_0^{t} \frac{ds}{g(s)^2} \leq   t +C_2
\la{g:lt}
\ee
for some $C_2>0$ and large $t>0$. Combining (\ref{def:Y}), (\ref{h:lt}) and (\ref{g:lt}) we conclude that (\ref{tn}) is valid provided $t>0$ is large enough.

Suppose now that $t>0$ is large enough so that (\ref{tn}) is valid. We then have
\[
Y(t) = (n-2)\frac{h(\rho)^2 \sinh t }{g(t)^2  \sinh \rho }  \geq (n-2)\frac{\sinh t }{ \sinh \rho } 
  \geq (n-2)\frac{\sinh t }{ \sinh (t+c_n) }  \longrightarrow (n-2) e^{-c_n} \; , \qquad \mbox{ as }t\to +\infty.
\]
We thus conclude that if $\alpha_n>0$ is such that $X(\alpha_n)< (n-2)e^{-c_n}$ then
\[
Y(t)  \geq X \big(  \alpha_n\tanh\frac{t}{2}  \big)  
\]
provided $t$ is large enough.

{\em Case 2: Small $t>0$.} It is not difficult to see that the constant $A$ in (\ref{A}) is negative. It then easily follows that
from (i) of Lemma \ref{lem:asympt_gh} that there exists $\mu>0$ so that
\be
\int_0^t \frac{ds}{g(s)^2}  \leq \frac{\pi^2}{ 2( -\ln 2t -\mu  )}  \; , 
\la{eq:int_g}
\ee
for small enough $t>0$.
We now distinguish cases according to the dimension $n$.

(i) $n\geq 5$. Applying Lemma \ref{lem:asympt_gh} we easily obtain
\be
\int_0^{\rho} \frac{ds}{h(s)^2}  = \frac{n-2}{2(c_2^{\#})^2} (2\rho)^{n-2} \big(    1-\tau_n \rho^2 +O(\rho^3)  \big)  \; , \qquad \mbox{ as }\rho\to  0+ \, ,
\la{eq:int_h}
\ee
where
\[
\tau_n =\frac{(n-1)(n-2)(n-3)}{12n(n-4)}.
\]
From (\ref{intro:fnr}),   (\ref{eq:int_g}) and (\ref{eq:int_h}) we obtain that for some $c>0$ and all small enough $t>0$  there holds,
\be
 (n-2)   \Big(\frac{ 1}{\pi c_2^{\#}}\Big)^2   (2\rho)^{n-2} \big(    1-\tau_n \rho^2 -c\rho^3 \big)   \leq  \frac{1}{ -\ln (2t)-\mu}.
\la{ex1}
\ee
The required inequality (\ref{yx}) can be written as
\be
\frac{ \sinh t}{g(t)^2} \geq  \frac{1}{n-2} X\big(  \alpha\tanh\frac{t}{2}  \big)   \frac{ \sinh \rho}{h(\rho)^2} \;\; \qquad \mbox{(small $t>0$).}
\la{ri}
\ee
We note that for small $t>0$ we have
\[
X\big(  \alpha\tanh\frac{t}{2}  \big)   \leq X\big( \frac{\alpha t}{2}  \big) = \frac{1}{ 1 - \ln(\alpha/4) -\ln (2t)}.
\]
Writing $M= 1 - \ln(\alpha/4)$ and using Lemma \ref{lem:asympt_gh}
we conclude that (\ref{ri}) will follow if we establish that for $M>0$ large enough there holds
\[
\frac{ \pi^2 }{ \big(\ln (2t))^2 +c_1 \ln (2t) \big)  } \geq \frac{1}{(n-2) (M -\ln (2t))} \cdot  \frac{ 1+\frac{\rho^2}{6} +c\rho^4} { \frac{ (c_2^{\#})^2}{ (n-2)^2} (2\rho)^{-n+2}\Big(  1+ \frac{(n-1)(n-3)}{12(n-4)}\rho^2 -c\rho^3  \Big)}
\]
where $c>0$ and  $c_1\in\R$ are fixed constants and the inequality is required to be valid for small enough $t>0$.
This is also written as
\be
(n-2) \Big(\frac{ 1}{ \pi c_2^{\#}}\Big)^2
 (2\rho)^{n-2} \cdot  \frac{1}{M-\ln (2t)} \leq  \frac{1}{   \big(  \ln^2 (2t) +c_1 \ln (2t) \big)  } \cdot \frac{ 1+ \frac{(n-1)(n-3)}{12(n-4)}\rho^2 -c\rho^4}{ 1+\frac{\rho^2}{6} +c\rho^4}.
\la{ex3}
\ee
From (\ref{ex1}) and (\ref{ex3}) we conclude that it is enough to establish the inequality
\[
\frac{  \ln^2 (2t) +c_1 \ln (2t)  }{ ( -\ln (2t) -\mu ) (M -\ln (2t)) } \leq 
\frac{ \big(1+ \frac{(n-1)(n-3)}{12(n-4)}\rho^2 -c\rho^3\big)\big(    1-\tau_n \rho^2 -c\rho^3  \big)}{ 1+\frac{\rho^2}{6} +c\rho^4}.
\]
Now, since $n\geq 5$,
\[
 \frac{1}{6} +\tau_n < \frac{(n-1)(n-3)}{12(n-4)} ,
\]
hence for small $\rho>0$,
\[
  \frac{ \Big(  1+ \frac{(n-1)(n-3)}{12(n-4)}\rho^2 -c\rho^3  \Big) \big(    1-\tau_n \rho^2  -c\rho^3 \big) }
  { 1+\frac{\rho^2}{6} +c\rho^4  } \geq 1.
\]
On the other hand it is easily seen that if $M>0$ is large enough then
\be
\frac{  \ln^2 (2t) +c_1 \ln (2t)  }{ ( -\ln (2t) -\mu ) (M -\ln (2t)) }
< 1 \,  ,
\la{29}
\ee
for small enough $t>0$. This completes the proof.

(ii) $n=4$. Applying Lemma \ref{lem:asympt_gh} we easily obtain
\[
\int_0^{\rho} \frac{ds}{h(s)^2}  = \frac{4}{(c_2^{\#})^2} \rho^2  \big(    1 +\frac{1}{8} \rho^2 \ln(2\rho) +O(\rho^2)  \big)  \; , \qquad \mbox{ as }\rho\to  0+ \, .
\]
Hence, by (\ref{eq:int_g}), for some $c>0$ and small $t>0$ there holds,
\be
   \Big(\frac{ 1}{\pi  c_2^{\#}}\Big)^2  (2\rho)^2 \big(    1 +\frac{1}{8} \rho^2 \ln(2\rho) -c\rho^2  \big)   \leq  \frac{1}{ 2\big( -\ln (2t)-\mu \big)}.
\la{ex1_4}
\ee
Arguing as in the case $n\geq 5$ we conclude that it is enough to establish that for $M>0$ large enough there holds
\[
2 \Big(\frac{ 1}{ \pi  c_2^{\#}}\Big)^2
 (2\rho)^{2} \cdot  \frac{1}{M-\ln (2t)} \leq  \frac{1}{   \big(  \ln^2 (2t) +c_1 \ln (2t) \big)  } \cdot \frac{ 1 - \frac{1}{4}\rho^2  
 \ln(2\rho)-c\rho^2}{ 1+c\rho^2},
\]
where $c>0$ is some fixed constant and the inequality is required to be valid for small enough $t>0$. 
Combining this with (\ref{ex1_4}) we conclude that it is enough to establish that for small $t>0$ there holds
\[
\frac{  \big(\ln (2t))^2 +c_1 \ln (2t) \big)}{ (\mu -\ln (2t)) (M -\ln (2t)) } \leq  
  \frac{ (1 +\frac{1}{8}\rho^2 \ln(2\rho)  -c \rho^2) ( 1 -\frac{1}{4}\rho^2 \ln(2\rho)  -c \rho^2)}{ 1 +c\rho^2 }.
\]
Since for small $\rho>0$ we have
\[
\frac{ (1 +\frac{1}{8}\rho^2 \ln(2\rho)  -c \rho^2) ( 1 -\frac{1}{4}\rho^2 \ln(2\rho)  -c \rho^2)}{ 1 +c\rho^2 } \geq 1,
\]
the result follows from (\ref{29}).

(iii) $n=3$. In this case we take $h(\rho)=1$ and therefore  the LHS of (\ref{eq:int_h}) is equal to $\rho$. The rest of the argument is similar, indeed, simpler, than that of the cases $n\geq 5$ or $n=4$.
This completes the proof of the theorem. $\hfill\Box$

{\bf Remark.} The function $X\big(\alpha \tanh(\rho/2) \big)$ captures the actual small time behaviour of $Y(t)$ in the sense that
\[
 \lim_{t\to 0+} \big( -Y(t) \ln t\big) =1 ;
\]
we omit the proof of this statement.

\subsection{Identifying the constant $\overline{S}_{3,p}$}

Our aim in this subsection is to prove that $\overline{S}_{3,p}=S_{3,p}$ (cf. (\ref{identify})) and give the proof of Theorem \ref{thm:abcde3}.
For this we shall use the half-space model of $\Hp^n$.

We start by establishing an inequality which is a consequence of inequality (\ref{ws}) and which will be used later on.
\begin{theorem}
Let $n\geq 3$, $2< p\leq 2^*$. Then for all  $u\in \cic(\R^n)$ there holds
\be
 \int_{\R^n}|\nabla u|^2dx   \geq
 S_{n,p}\bigg( \int_{\R^n}  \Big(  \frac{ |x-e_n| \; |x+e_n|}{2} \Big)^{ \frac{p(n-2)}{2}-n} |u|^{p}  dx \bigg)^{2/p} .
\la{501}
\ee
Moreover the constant is sharp and is attained by the function
\[
u(x) =\Big(  |x+e_n|^{\frac{(p-2)(n-2)}{2}}  +  |x-e_n|^{\frac{(p-2)(n-2)}{2}} \Big)^{-\frac{2}{p-2}} , \qquad x\in\R^n\, .
\]
\la{thm:ten}
\end{theorem}
{\em Proof.} The map
\[
S(y)= \frac{1}{|y +e_n|^2} (2y',1-|y|^2)
\]
maps conformally $\R^n$ onto $\R^n$. We note that
\[
|S(y)| =\frac{|y-e_n|}{|y+e_n|}.
\]
The Jaccobian determinant $JS(y)$ of $S$ can be computed explicitly and one finds
\[
|JS(y)| =\frac{2^n}{ |y +e_n|^{2n}} .
\] 
Now let $u\in C^{\infty}_c(\R^n)$ be given. We define the function $w$ by
\[
w(y) =u(S(y))  | (JS)(y)|^{\frac{n-2}{2n}}  = u(S(y)) \Big( \frac{2}{|y+e_n|^2} \Big)^{\frac{n-2}{2}}.
\]
By (\ref{ws}) we have
\be
\int_{\R^n}|\nabla w|^2 dy \geq S_{n,p} \Big(  \int_{\R^n}|y|^{ \frac{p(n-2)}{2}-n}|w|^p dy\Big)^{2/p} .
\la{ws1}
\ee
Changing variables via $S$, $x=S(y)$, in (\ref{ws1}) we arrive at (\ref{501}).

Finally, under this transformation the minimizer (\ref{minim}) is transformed to the function
\bean
w(y) &=& |JS(y)|^{\frac{n-2}{2n}} u(x)   \\
&=&   \frac{2^\frac{n-2}{2}}{ |y +e_n|^{n-2}}  \big(  1+|x|^{\frac{(p-2)(n-2)}{2}} \big)^{-\frac{2}{p-2}} \\
&=& 2^{\frac{n-2}{2}} \Big(  |y+e_n|^{\frac{(p-2)(n-2)}{2}}  +  |y-e_n|^{\frac{(p-2)(n-2)}{2}} \Big)^{-\frac{2}{p-2}}.
\eean
This concludes the proof. $\hfill\Box$

In case $n=3$ we have
\begin{theorem}
For all $2<p\leq 6$ there holds
\[
\overline{S}_{3,p} =S_{3,p} = \frac{p}{2^{\frac{2}{p}}} 
\bigg[ 
\frac{ 4 \pi \Gamma^2(\frac{p}{p-2})}{ (p-2)  \Gamma(\frac{2p}{p-2}) } 
 \bigg]^{\frac{p-2}{p}}.
\]
\la{thm15}
\end{theorem}
{\em Proof.} We begin by recalling from \cite{BFL} that when $n=3$ the fundamental solution of the equation
\be
u_t =\Delta u  + \frac{1}{4x_3^2} \; , \qquad x=(x',x_3)\in\R^3_+ \; , \;\; t>0,
\la{semi}
\ee
is given by
\[
G(x'-y' ,x_3,y_3,t)=  \frac{1}{(4\pi t)^2} \sqrt{x_3y_3} \;  e^{ -\Frac{ |x'-y'|^2 +x_3^2 +y_3^2}{4t}}
\int_{0}^{2\pi} e^{ \Frac{x_3y_3}{2t} \cos\phi } \, d\phi \; .
\]
Let $Q$ denote the generator of the semigroup associated to (\ref{semi}) and $Q^{-1}(x,y)$ be the integral kernel of $Q^{-1}$ so that.
\[
Q^{-1}(x,y) = \int_0^{\infty} G(x'-y' ,x_3,y_3,t) dt 
\]
It has been proved in \cite{BFL} that we then have the estimate
\[
Q^{-1}(x,y) \leq  (-\Delta)^{-1}(x,y)  \; , \quad\quad x,y\in\R^3_+,
\]
where
\[
 (-\Delta)^{-1}(x,y) =  \frac{1}{4\pi |x-y|} 
\]
is the Green function for $-\Delta$ in $\R^3$.

Now, from Theorem \ref{thm:ten} we have
\[
 \int_{\R^3}|\nabla u|^2dx   \geq
 S_{3,p}\bigg( \int_{\R^3}  \Big(  \frac{ |x-e_3| \; |x+e_3|}{2} \Big)^{ \frac{p-6}{2}} |u|^{p}  dx \bigg)^{2/p} .
\]
This is also written as
\be
\inprod{-\Delta u}{u} \geq S_{3,p} \|  \phi^{-1}u  \|_p^2
\la{505}
\ee
where
\[
\phi(x) = \Big(  \frac{ |x-e_3| \; |x+e_3|}{2} \Big)^{\frac{6-p}{2p}} \; .
\]
Using duality (\ref{505}) gives
\[
S_{3,p}   \inprod{(-\Delta)^{-1} v}{v} \leq \|  \phi v  \|_{p'}^2 \; .
\]
We conclude as in \cite{BFL} that
\bean
\inprod{ \phi^{-1}u}{\phi v}^2&=& \inprod{ u}{ v}^2 = \inprod{Q^{1/2}u}{Q^{-1/2}v}^2 \leq \inprod{Qu}{u}  \inprod{Q^{-1}v}{v} \\
&\leq & \inprod{Qu}{u}  \inprod{ (-\Delta)^{-1}v}{v} \leq \frac{1}{S_{3,p}} \inprod{Qu}{u}  \|  \phi v  \|_{p'}^2 \, .
\eean
Therefore
\[
 \| \phi^{-1}u\|_p^2 \leq \frac{1}{S_{3,p}} \inprod{Qu}{u} ,
\]
that is
\[
\int_{\R^3_+}|\nabla u|^2 dx -\frac{1}{4}\int_{\R^3_+}\frac{u^2}{x_3^2}dx \geq S_{3,p}\bigg( \int_{\R^3_+}  \Big(  \frac{ |x-e_3| \; |x+e_3|}{2} \Big)^{ \frac{p-6}{2}} |u|^{p}  dx \bigg)^{2/p}.
\]
Hence $\overline{S}_{3,p} \geq S_{3,p}$.
The reverse inequality $\overline{S}_{3,p} \leq S_{3,p}$ follows by noting that $S_{3,p}$ is the best constant for the inequality
\[
 \int_{\R^3_+}|\nabla u|^2dx  \geq 
  S_{3,p}\bigg( \int_{\R^3_+}  \Big(  \frac{ |x-e_3| \; |x+e_3|}{2} \Big)^{ \frac{p-6}{2}} |u|^{p}  dx \bigg)^{2/p} , \quad\quad u\in \cic(\R^3_+).
\]
$\hfill\Box$

\

{\bf \em Proof of Theorem \ref{thm:abcde3}.} Inequality (\ref{ell}) follows from Theorems \ref{thm:p_h_s} and \ref{thm15}.
The sharpness of $S_{3,p}$ follows by a local argument near the origin. One uses the fact that the Sobolev constant of Theorem \ref{thm:2019} remains invariant if we restrict to test functions with support in any given small neighbourhood of the origin.  $\hfill\Box$

\section{Hardy-Sobolev inequalities with a boundary point singularity}
\la{section_bps}

In this section  we obtain Hardy-Sobolev inequalities when we place a point singularity on the boundary of a bounded domain $\Omega\subset\R^n$. Before doing so we consider the flat case $\Omega=B_1^+=\{x\in B_1 \, : \; x_n>0\}$.

Given $n\geq 3$, $0\leq \gamma< n/2$ and  $2<p\leq 2^*$ we define
\be
S^*_{n,p,\gamma} =\inf_{\cic(\R^n_+)}\frac{
\Int_{\R^n_+}|\nabla u|^2dx - \gamma(n-\gamma)\Int_{\R^n_+}\frac{u^2}{|x|^2}dx }
{ \Big( \Int_{\R^n_+}  |x|^{\frac{p(n-2)}{2}-n} |u|^{p}dx\Big)^{2/p}} .
\la{kal1}
\ee
Our first result reads
\begin{theorem}
Let $n\geq 3$, $2< p\leq 2^*$ and $0\leq \gamma <n/2$ be given. Let
\[
\alpha_{n,\gamma}= e^{\frac{n-1-2\gamma}{n-2\gamma}}.
\]
Then for all $0<\alpha\leq \alpha_{n,\gamma}$ and all $u\in \cic(B_1^+)$ there holds
\be
\int_{B_1^+}|\nabla u|^2 dx  -\frac{n^2}{4} \int_{B_1^+}\frac{u^2}{|x|^2} dx
 \geq  (n-2\gamma)^{-\frac{p+2}{p}}S_{n,p,\gamma}^* \bigg(  \int_{B_1^+}  |x|^{\frac{p(n -2)}{2} -n} X^{\frac{p+2}{2}}(\alpha |x|) |u|^{p} dx  \bigg)^{2/p}.
\la{8881}
\ee
\end{theorem}
{\em Proof.} 
 A simple scaling argument shows that the best constant $S^*_{n,p,\gamma}$ in (\ref{kal1}) remains the same if $\R^n_+$ is replaced by $B_{\rho}^+$ for any $\rho>0$.
Making the change of variables
\[
u(x) = \frac{x_n}{|x|^{\gamma}} v(x)
\]
we then obtain
\bean
S^*_{n,p,\gamma} &=& \inf\frac{\Int_{B_{\rho}^+ } \frac{x_n^2}{|x|^{2\gamma}} |\nabla v|^2dx}
{ \Big( \Int_{B_{\rho}^+}   x_n^p   |x|^{ \frac{ p(n-2-2\gamma)}{2} -n}  |v|^{p}dx\Big)^{2/p}}  \\
&=& \inf \frac{ \Int_0^{\rho}\Int_{\sphp}\omega_n^2  r^{n+1-2\gamma}   \Big( v_r^2  +
\frac{1}{ r^2} |\nabla_{\omega} v|^2  \Big)dS(\omega) \, dr  } {
\Big(  \Int_0^{\rho}\Int_{\sphp}   \omega_n^{p}  r^{ \frac{p(n-2\gamma)}{2}-1 }  |v|^{p}  dS(\omega) \, dr  \Big)^{2/p }}  ;  
\eean
here $\omega_n$ denotes the $n$th component of $\omega\in \sphp$.
We next change variables setting
\[
t= \frac{1}{ n-2\gamma} r^{2\gamma -n} \; , \qquad v(r,\omega) =w(t,\omega).
\]
After some more computations we arrive at
\be
(n-2\gamma)^{- \frac{p+2}{p}}  S^*_{n,p,\gamma}  =
\inf\frac{\Int_{\frac{\rho^{2\gamma-n}}{n-2\gamma}}^{\infty}\Int_{\sphp}   \omega_n^2 \Big( w_t^2  + 
\frac{1}{ (n-2\gamma)^2t^2 } |\nabla_{\omega} w|^2  \Big)dS(\omega) \, dt  }
{\Big(  \Int_ {\frac{\rho^{2\gamma-n}}{n-2\gamma}} ^{\infty}\Int_{\sphp}  \omega_n^{p} \, t^{-\frac{p+2}{2}}   |w|^{p}  dS(\omega) \, dt  \Big)^{2/p }} .
\la{kal5}
\ee
It is enough to establish (\ref{8881}) for $\alpha=\alpha_{n,\gamma}$.
Let
\[
\tau^*_{n,p,\gamma} =\inf
\frac{\Int_{B_1^+}|\nabla u|^2 dx  -\frac{n^2}{4} \Int_{B_1^+}\frac{u^2}{|x|^2} dx}
{\bigg(  \Int_{B_1^+}  |x|^{\frac{p(n -2)}{2} -n} X^{\frac{p+2}{2}}(\alpha |x|) |u|^{p} dx  \bigg)^{2/p}}
\]
Setting
\[
u(x) =\frac{x_n}{|x|^{\frac{n}{2}}} v(x)
\]
we find
\[
\tau^*_{n,p,\gamma} 
= \inf\frac{\Int_0^1 \Int_{\sphp}  r\omega_n^2  \Big(  v_r^2  + \frac{1}{r^2} |\nabla_{\omega}v|^2 \Big) dS \, dr}
 { \Big( \Int_0^1 \Int_{\sphp} \omega_n^{p} \, r^{-1}   X^{\frac{p+2}{2}}(\alpha r) |v|^{p} dS \, dr \Big)^{2/p}}.
\]
We change variables again setting
\[
t= \frac{1}{X(\alpha r)} \; , \qquad v(r,\omega) =w(t,\omega).
\]
After some more computations we obtain
\be
\tau^*_{n,p,\gamma} =\inf \frac{ \Int_{\frac{1}{X(\alpha)}}^{+\infty} \Int_{\sphp}  \omega_n^2  \big(  w_t^2  +|\nabla_{\omega}w|^2 \big) dS \, dt}
{ \Big( \Int_{\frac{1}{X(\alpha)}}^{+\infty} \Int_{\sphp} \omega_n^p \;   t^{-\frac{p+2}{2}} |w|^{p} dS \, dt \Big)^{2/p}}.
\la{1234561}
\ee
Choosing $\rho=1$ and noting that $X(\alpha_{n,\gamma}) =n-2\gamma$, we compare (\ref{kal5}) and (\ref{1234561}) and obtain $\tau^*_{n,p,\gamma}
\geq  (n-2\gamma)^{-\frac{p+2}{p}}S_{n,p,\gamma}^*$ as required. $\hfill\Box$


Actually inequality (\ref{8881}) can be improved. The next result plays an important role in establishing Theorem \ref{c} which is the main result of this section.
\begin{theorem}
Let $n\geq 3$, $2< p\leq 2^*$, $0\leq \gamma < n/2$ and $0<\theta <2$. Let $R=R_{\theta,\gamma}$ and $\alpha_{n,\gamma,\theta}$ be defined by
\[
R^{\theta} =1 +\frac{1}{\sqrt{n-2\gamma}}
  \;\; , \qquad  -\ln\alpha_{n,\gamma, \theta}= 
R^{2\theta} -1 + \int_0^1 \frac{s^{\theta-1}(  2R^{\theta} -s^{\theta} )}{ (R^{\theta} -s^{\theta})^2}ds
\]
Then for all $0<\alpha\leq \alpha_{n,\gamma, \theta}$ and for all $u\in \cic(B_1^+)$ there holds
\begin{eqnarray}
&&\hspace{-3cm}\int_{B_1^+}|\nabla u|^2 dx  -\frac{n^2}{4} \int_{B_1^+}\frac{u^2}{|x|^2} dx
 -\theta^2 \int_{B_1^+} \frac{ u^2}{ |x|^{2-\theta} (R^{\theta} -|x|^{\theta})} dx  \nonumber \\[0.2cm]
&\geq&  (n-2\gamma)^{-\frac{p+2}{p}}S_{n,p,\gamma}^* \bigg(  \int_{B_1^+}  |x|^{\frac{p(n -2)}{2} -n} X^{\frac{p+2}{2}}(\alpha |x|) |u|^{p} dx  \bigg)^{2/p}.
\la{888}
\end{eqnarray}
\la{thm:el}
\end{theorem}
{\em Proof.} We recall (cf. (\ref{beta})) that
\[
B(r)  =\frac{1}{ (  R^{\theta} -r^{\theta} )^2 \Big(  1+ \int_r^1 \frac{dt}{t(  R^{\theta} -t^{\theta} )^2} \Big)} \; , \qquad r\in (0,1).
\]
To prove (\ref{888}) we shall first establish the following inequality for all $u\in \cic(B_1^+)$
\bea
&& \int_{B_1^+}|\nabla u|^2 dx  -  \frac{n^2}{4} \int_{B_1^+}\frac{u^2}{|x|^2} dx -\theta^2 \int_{B_1^+} \frac{ u^2}{ |x|^{2-\theta} 
(R^{\theta} -|x|^{\theta})} dx   \nonumber \\
&& \hspace{3cm} \geq  (n-2\gamma)^{- \frac{p+2}{p}} S_{n,p,\gamma}^* \bigg(  \int_{B_1^+} |x|^{ \frac{p(n-2)}{2}-n} B^{\frac{p+2}{2}}(|x|) |u|^{p}dx  \bigg)^{2/p}.
\la{pinbs00}
\eea
So let us define
\[
\tau^*_{n,p,\gamma,\theta} =\inf \frac{\Int_{B_1^+}|\nabla u|^2 dx  -  \frac{n^2}{4} \Int_{B_1^+}\frac{u^2}{|x|^2} dx -\theta^2 \Int_{B_1^+} \frac{ u^2}{ |x|^{2-\theta} 
(R^{\theta} -|x|^{\theta})} dx}{ \bigg(  \Int_{B_1^+} |x|^{ \frac{p(n-2)}{2}-n} B^{\frac{p+2}{2}}(|x|) |u|^p dx  \bigg)^{2/p}}.
\]
Setting
\[
u(x) =\frac{x_n}{|x|^{\frac{n}{2}}} (R^{\theta} -|x|^{\theta}) v(x)
\]
we find
\bean
\tau^*_{n,p,\gamma,\theta} &=&\inf\frac{ \Int_{B_1^+}  \frac{x_n^2}{|x|^n} (R^{\theta} -|x|^{\theta})^2 |\nabla v|^2 dx}{\bigg( \Int_{B_1^+} 
x_n^{p}|x|^{-p-n} (R^{\theta} -r^{\theta})^{p}  B^{\frac{p+2}{2}}(|x|) |v|^{p} dS \, dr \bigg)^{2/p}}   \\
&=& \inf\frac{\Int_0^1 \Int_{\sphp}  r\omega_n^2  (R^{\theta} -r^{\theta})^2 \Big(  v_r^2  + \frac{1}{r^2} |\nabla_{\omega}v|^2 \Big) dS \, dr}
 { \Big( \Int_0^1 \Int_{\sphp} \omega_n^{p} \, r^{-1} (R^{\theta} -r^{\theta})^{p}  B^{\frac{p+2}{2}}(r) |v|^{p} dS \, dr \Big)^{2/p}}.
 \eean
We next change variables by
\be
 t= 1+ \int_r^1 \frac{ds}{s(  R^{\theta} -s^{\theta} )^2}             \;\; , \qquad   v(r,\omega) =w(t,\omega)
\la{t:r}
\ee
and obtain
\be
\tau^*_{n,p,\gamma,\theta} =\inf \frac{ \Int_1^{+\infty} \Int_{\sphp}  \omega_n^2  \big(  w_t^2  +Q(t) |\nabla_{\omega}w|^2 \big) dS \, dt}
{ \Big( \Int_1^{+\infty} \Int_{\sphp} \omega_n^p \;   t^{-\frac{p+2}{2}} |w|^{p} dS \, dt \Big)^{2/p}}
\la{123456}
\ee
where
\[
Q(t)=(R^{\theta} -r^{\theta})^4
\]
and $r=r(t)$ is the inverse of (\ref{t:r}).
Choosing $\rho$ in (\ref{kal5}) so that $\rho^{n-2\gamma} =n-2\gamma$ we obtain
\be
(n-2\gamma)^{- \frac{p+2}{p}}  S^*_{n,p,\gamma}  =
\inf\frac{\Int_{1}^{\infty}\Int_{\sphp}   \omega_n^2 \Big( w_t^2  + 
\frac{1}{ (n-2\gamma)^2t^2 } |\nabla_{\omega} w|^2  \Big)dS(\omega) \, dt  }
{\Big(  \Int_ {1} ^{\infty}\Int_{\sphp}  \omega_n^{p} \, t^{-\frac{p+2}{2}}   |w|^{p}  dS(\omega) \, dt  \Big)^{2/p }} .
\la{kal555}
\ee

By the choice of $R$ we have
\[
 Q(t)\geq  \frac{1}{(n-2\gamma)^2t^2}  \; , \qquad t\geq 1. 
\] 
Therefore from (\ref{123456}) and (\ref{kal555}) we obtain that $\tau^*_{n,p,\gamma,\theta} \geq (n-2\gamma)^{-\frac{p+2}{p}}S_{n,p,\gamma}^*$ and
(\ref{pinbs00}) follows. 
Finally, inequality (\ref{888}) now follows by recalling from Lemma \ref{lem:bx} that $B(r)\geq X(\alpha_{n,\gamma,\theta}\; r)$ for $r\in (0,1)$.
$\hfill\Box$

\

We can now state and prove the main result of this section. In the particular case $p=2^*$ this is Theorem \ref{c}.
\begin{theorem}
\label{c_p}
Let $\Omega\subset\R^n$, $n\geq 3$, be a bounded domain with $0 \in \partial \xO$  and let $D=\sup_{\Omega}|x|$.
Assume that $\Omega$ satisfies an exterior ball condition at zero with exterior ball $B_{\rho}(-\rho e_n) \subset \cC \overline{\xO} $. 
Then for any $2<p\leq 2^*$ and any $\gamma\in [0,n/2)$ there exist an $r_{n,\gamma} $ and $\alpha_{n,\gamma}^*$  in $(0,1)$  such that,
if the radius $\rho$  of the exterior ball satisfies  $\rho\geq D/ r_{n,\gamma}$ then for all $0<\alpha\leq\alpha_{n,\gamma}^*$ there holds
\be
\int_{\xO} |\nabla u|^2 dx \geq   \frac{n^2}{4}  \int_{\xO} \frac{u^2}{|x|^2} dx \\
+  (n-2\gamma)^{- \frac{p+2}{p}}  S^*_{n,p,\gamma} \bigg( \int_{\xO} |x|^{\frac{p(n-2)}{2}-n} \Big(  \frac{|x+2\rho e_n|}{2\rho} \Big)^{\frac{p(n-2)}{2}-n}   X^{\frac{p+2}{2}}|u|^{p}     dx \bigg)^{\frac{2}{p}},
\la{c111}
\ee
for all $u \in C^{\infty}_{c}(\xO)$; here $X=X(\alpha |x|/ D)$.
\end{theorem}
{\em Proof.} We shall establish the result in case $\rho=1$, the general case will then follow by scaling. For simplicity of subsequent computations we  make a translation by $e_n$ and place the singularity at $e_n$ so that the exterior ball is $B_1(0)$; hence instead of  (\ref{c111}) we shall establish
\[
\int_{\xO} |\nabla u|^2 dx \geq   \frac{n^2}{4}  \int_{\xO} \frac{u^2}{|x-e_n|^2} dx \\
+  (n-2\gamma)^{- \frac{p+2}{p}}  S^*_{n,p,\gamma} \bigg( \int_{\xO} |x-e_n|^{\frac{p(n-2)}{2}-n} \Big(  \frac{|x+ e_n|}{2} \Big)^{\frac{p(n-2)}{2}-n}   X^{\frac{p+2}{2}}|u|^{p}     dx \bigg)^{\frac{2}{p}}.
\]

As a first step we shall establish that there exist $r =r_{n,\gamma} \in (0,1)$  such that
\bea
&& \int_{\cC B_1 \cap B( e_n,r  )}|\nabla u|^2dx \geq \frac{n^2}{4}\int_{\cC B_1 \cap B( e_n, r)}
\frac{u^2}{|x- e_n|^2}dx   \nonumber \\
&& \hspace{.5cm} +   (n-2\gamma)^{- \frac{p+2}{p}}  S^*_{n,p,\gamma} \bigg( \int_{\cC B_1 \cap B( e_n,r  )} |x-e_n|^{\frac{p(n-2)}{2}-n} \Big(  \frac{|x+ e_n|}{2} \Big)^{\frac{p(n-2)}{2}-n}   X^{\frac{p+2}{2}}|u|^{p} dx \bigg)^{\frac{2}{p}} ,
\la{gera:eq:1.33aaa} 
\eea
for all $u\in C^{\infty}_c(\cC B_1 \cap B( e_n, r ))$; here $X=X( \alpha_{n,\gamma}^* |x-e_n|)$, where
 $\alpha_{n,\gamma}^*$ is a constant that depends only on $n$ and $\gamma$.


To prove this we consider the conformal map
\[
T(y) =  \frac{1}{ |y-e_n|^2}(2y',1-|y|^2) 
\]
which maps $\R^n_+$ onto $\cC B_1$.
The Jacobian of $T$ is easily computed and we find
\be
|JT(y)| = \frac{2^n}{|y-e_n|^{2n}}.
\la{gera:jt1}
\ee
In addition we have
\[
T^{-1}(x)=\frac{1}{|x'|^2+(x_n+1)^2}(2x',|x|^2 -1)
\]  
and therefore
\be
|T^{-1}(x)| =\frac{|x-e_n|}{|x+e_n|}.
\la{gera:T1112}
\ee
Now let $r<1$ be fixed (this will be chosen later on)
and let $F\in C^{\infty}_c(T(B_r^+))$ be given. We define the function $G$ on $B_r^+$ by
\[
G(y) =F(T(y))  | (JT)(y)|^{\frac{n-2}{2n}}  = F(T(y)) \Big( \frac{2}{|y-e_n|^2} \Big) ^{\frac{n-2}{2}}.
\]
Applying Theorem \ref{thm:el} with $\theta=1/2$ (hence $\sqrt{R} =1 +\frac{1}{\sqrt{n-2\gamma}}$) and
using a scaling argument we obtain that there exists $\alpha_{n,\gamma}$ such that
\bea
&& \int_{B_r^+}|\nabla G|^2  dy \geq \frac{n^2}{4} \int_{B_r^+}\frac{G^2}{|y|^2} dy  +
 \frac{1}{4R^{1/2}r^{1/2}}\int_{B_{r}^+}\frac{G^2}{|y|^{3/2}}dy   \nonumber \\
&& \hspace{2.5cm} \geq (n-2\gamma)^{- \frac{p+2}{p}}  S^*_{n,p,\gamma}  
\bigg(  \int_{B_r^+}  |y|^{\frac{p(n -2)}{2} -n} X^{\frac{p+2}{2}} |G|^{p} dy  \bigg)^{2/p}  ,
\label{gera:josbir1}
\eea
where $X=X(\alpha_{n,\gamma} |y|/r)$.

We have (cf. \cite{BFT})
\[
\int_{B_r^+}|\nabla G|^2  dy  =  \int_{T(B_r^+)}|\nabla F|^2  dx \, .
\]
Using  (\ref{gera:jt1}) and (\ref{gera:T1112}) we also find that
\[
\int_{B_r^+}\frac{G^2}{|y|^2} dy = \int_{T(B_r^+)}\frac{4F^2}{|x-e_n|^2|x+e_n|^2} dx \, .
\]
The other two integrals in (\ref{gera:josbir1}) can similarly be transformed and we conclude that (\ref{gera:josbir1})
takes the following equivalent form
\bea
&& \int_{T(B_r^+)}|\nabla F|^2dx \geq \frac{n^2}{4} \int_{T(B_r^+)}\frac{4F^2}{|x-e_n|^2|x+e_n|^2}dx   + \frac{1}{4R^{1/2}r^{1/2}} \int_{T(B_r^+)}\frac{4F^2}{|x-e_n|^{3/2} |x+e_n|^{5/2}}dx  \nonumber \\
&& \hspace{1cm}  +  (n-2\gamma)^{- \frac{p+2}{p}}  S^*_{n,p,\gamma} \bigg( \int_{T(B_r^+)} |x-e_n|^{\frac{p(n-2)}{2}-n} \Big(  \frac{|x+ e_n|}{2} \Big)^{\frac{p(n-2)}{2}-n}   X^{\frac{p+2}{2}}|F|^{p} dx \bigg)^{\frac{2}{p}}
\la{ipp}
\eea
where $X=X(\alpha_{n,\gamma} |x-e_n| /r|x+e_n|)$.
Now, it follows from (\ref{gera:T1112}) and some simple geometry that for any $r<1$
\[
T(B_r)=\{ x\in\R^n \, : \, |x-e_n| <  r|x+e_n|\} =B_{\frac{2r}{1-r^2}}\big(\frac{1+r^2}{1-r^2} e_n \big) \supset B_r(e_n),
\]
therefore
\be
T(B_r^+)  \supset B_1^{c}  \cap B_r(e_n)  .
\la{inclusion}
\ee
We will choose $r \in (0,1)$ such that for all $x\in B_1^{c}\cap B_r(e_n)
\subset T(B_r^+)$ there holds
\[
\frac{n^2}{4}\frac{4}{|x-e_n|^2|x+e_n|^2} +  \frac{1}{4R^{1/2}r^{1/2}}\frac{4}{|x-e_n|^{3/2} |x+e_n|^{5/2}}  
\geq \frac{n^2}{4|x-e_n|^2} ,
\]
or equivalently,
\be
|x-e_n|^{1/2} \geq  n^2R^{1/2}r^{1/2}|x+e_n|^{5/2} \Big( 1-\frac{4}{|x+e_n|^2}\Big).
\la{gef}
\ee
Indeed, this is immediate for $|x+e_n| \leq 2$.  Assuming that $|x+e_n| > 2$ we 
set $t=|x-e_n|$. We then have $|x+e_n|\leq t+2$ and therefore (\ref{gef}) will follow provided
\[
n^2R^{1/2} r^{1/2} t^{1/2}(t+4)(t+2)^{1/2} \leq 4,
\]
for all $t\leq r$. Simple computations give that the last inequality holds true provided $ t\leq 1 / (75n^4Rr)$.
This will be true for all $t\leq r$ if we choose
\[
r=r_{n,\gamma} :=\frac{1}{n\sqrt{75R}}.
\]
Finally, the inequality $|x+e_n| \leq 3$ implies $X(\alpha_{n,\gamma} |x-e_n| /r_{n,\gamma}|x+e_n|)\geq X( \alpha_{n,\gamma} |x-e_n|/3r_{n,\gamma})$.
Inequality (\ref{gera:eq:1.33aaa}) now follows with $\alpha_{n,\gamma}^* =  \alpha_{n,\gamma} /3r_{n,\gamma}$ by recalling (\ref{ipp}) and (\ref{inclusion}).

Since $D\leq r_{n,\gamma}$ we may choose $r=D$ in (\ref{gera:eq:1.33aaa}). Combining (\ref{gera:eq:1.33aaa}) with the inclusions
\[
\Omega\subset  \Omega \cap B(e_n,D) \subset \cC B( 1) 
\cap B(e_n ,D)   
\]
completes the proof of the theorem. $\hfill\Box$


\end{document}